\documentclass[12pt]{article}

\usepackage{amsfonts}
\usepackage{graphics}
\usepackage{amsmath}
\usepackage{amscd}
\usepackage{rlepsf}

\newtheorem{Theorem}{Theorem}

\newtheorem{Corollary}[Theorem]{Corollary}
\newtheorem{Example}[Theorem]{Example}
\newtheorem{Definition}[Theorem]{Definition}

\newtheorem{Remark}[Theorem]{Remark}

\newtheorem{Conjecture}{Conjecture}
\newtheorem{Lemma}[Theorem]{Lemma}
\newtheorem{Proposition}[Theorem]{Proposition}

\newtheorem{Exercise}{Exercise}

\newtheorem{Fundamental Theorem}{Fundamental Theorem}

\newenvironment{Proof}[1][Proof]{\textbf{#1.} }{\ \rule{0.5em}{0.5em}}

\def \l {{\leq t}}
\def \n {\nu}
\def \free {\mathrm{free}}
\def \Hom {\mathrm{Hom}}
\def \Z {\mathbb{Z}}
\def \id {\mathrm{id}}

\def \f {\phi}
\def \p {\psi}

\def \ra {\xrightarrow}
\def \cH {\mathcal{H}}
\def \H {\mathrm{Hom}}

\def \g {\gamma}
\def \G {\mathcal{G}}
\def \d {\partial}
\def \t {\triangleright}

\def \S {\Sigma}
\def \R {\mathbb{R}}
\def \s {\scriptstyle}

\def \N {\mathbb{N}_0}
\begin{document}

\title{On 2-Dimensional Homotopy Invariants of Complements of Knotted Surfaces}

\author{ Jo\~{a}o  Faria Martins\\ \footnotesize\it  {Departamento de Matem\'{a}tica, Instituto Superior T\'{e}cnico,}\\ {\footnotesize\it Av. Rovisco Pais, 1049-001 Lisboa, Portugal}\\ {\footnotesize\it jmartins@math.ist.utl.pt}}

\date{\today}

\maketitle
\abstract{
We prove that if $M$ is a CW-complex and $*$ is a 0-cell of $M$, then the crossed module $\Pi_2(M,M^1,*)=(\pi_1(M^1,*),\pi_2(M,M^1,*),\d,\t)$ does not depend on the cellular decomposition of $M$ up to free products with $\Pi_2(D^2,S^1,*)$, where $M^1$ is the 1-skeleton of $M$. From this it follows that  if $\G$ is a finite crossed module and  $M$ is finite, then the number of crossed module morphisms  $\Pi_2(M,M^1,*) \to \G$ (which is finite)  can be re-scaled to a homotopy invariant $I_\G(M)$ (i. e. not dependent on the cellular decomposition of $M$), a construction similar to David Yetter's in
  \cite{Y}, or Tim Porter's in \cite{P,P2}. We describe an
  algorithm to calculate $\pi_2(M,M^{(1)},*)$ as a crossed module over
  $\pi_1(M^{(1)},*)$, in the case when $M$ is the complement of a knotted surface
  in $S^4$ and $M^{(1)}$ is the $1$-handlebody of a handle decomposition of $M$, which, in particular, gives a method to calculate the algebraic 2-type of $M$. In addition, we prove that   the invariant $I_\G$ yields a non-trivial invariant of knotted surfaces. }

\tableofcontents
\section*{Introduction}
Let $(M,N,*)$ be a pair of based path-connected spaces.
The concept of a crossed module arises from  a universal description of the properties of the border map  $\d:\pi_2(M,N,*) \to \pi_1(N,*)$, together   with the natural action of $\pi_1(N,*)$ on $\pi_2(M,N,*)$. These  data define the crossed module $\Pi_2(M,N,*)$, called ``fundamental crossed module of $(M,N,*)$''.

 Due to some strong theorems by  J.H.C. Whitehead, it is possible, in principle, to calculate $\Pi_2(M,M^1,*)$, where   $M$ is a connected CW-complex, $M^1$ is its 1-skeleton and $*$ is a 0-cell of $M$.  The calculability of fundamental crossed modules is, in addition, enhanced by a 2-dimensional van Kampen theorem due to Ronald Brown and Philip Higgins. Note that the crossed module $\Pi_2(M,M^1,*)$ determines not only $\pi_1(M,*)$ and $\pi_2(M,*)$ as a module over $\pi_1(M,*)$,  but also it determines the k-invariant $k(M,*)\in H^3(\pi_1(M,*),\pi_2(M,*))$; in other words all the algebraic 2-type of $M$, thus being stronger than $\pi_1(M,*)$ and $\pi_2(M,*)$ alone.

Despite the asymmetry introduced by a choice of a 1-skeleton, the crossed module $\Pi_2(M,M^1,*)$ does not depend on the cellular decomposition of $M$ up to free products with  $\Pi_2(D^2,S^1,*)$. Therefore, in particular, if $\G$ is a finite crossed module and $M$ is a finite CW-complex then the number of crossed module morphisms from $\Pi_2(M,M^1,*)$ into $\G$ (which is finite) can be re-scaled to a homotopy invariant $I_\G(M)$.

We will use this set of results to develop an algorithm to calculate $\Pi_2(M,M^{(1)},*)$ in the case when $M$ is the complement of a knotted surface $\S$ in $S^4$ and $M^{(1)}$ is the 1-handlebody of a handle decomposition of it. Both the handle decomposition of $M$ (following \cite{G})  and the method  to determine $\Pi_2(M,M^{(1)},*)$ are defined  from a movie of $\S$.  We show how to calculate the second fundamental group of the Spun Trefoil complement from our construction.

In addition, we will prove that the counting invariant $I_\G$ defines a non trivial invariant of knotted surfaces, with good calculability, especially in the case of abelian crossed modules. This invariant has an innovative incorporation of information at saddle and maximal points of movies of knotted surfaces, which, should be able to be generalised for quandles, for example. In \cite{FM} we defined an invariant of knotted surfaces from any  finite crossed module. The
construction was inspired by previous work of David Yetter and Tim Porter on
manifold invariants defined from models of homotopy $2$-types (crossed modules
of groups), see \cite{Y} and \cite{P,P2}. This
article should give, in particular, an interpretation of our previous construction.
 \section{Preliminaries and General Results}

\subsection{Crossed Modules}
Let $G$ and $E$ be groups. A  crossed module with base $G$ and fibre $E$, say $\G=(G,E,\d,\t)$, is given by a group morphism $\d : E \to G$ and an action $\t$ of $G$ on $E$ on the left by automorphisms. The conditions on $\t$ and $\d$ are:
\begin{enumerate}
\item $\d(X \t e)=X\d(e)X^{-1};\forall X \in G, \forall e \in E$,
\item $\d(e) \t f=e f e^{-1}; \forall e, f \in E$.
\end{enumerate}
Notice that the second condition implies that $\ker \d$ commutes with all $E$.
We call $G$ the base group and $E$ the principal group. A crossed module is
called finite if both $G$ and $E$ are finite. A pre-crossed module is defined similarly to a crossed module, but skipping condition $2$.

The significance of the definition of crossed modules for Geometric Topology stems  from: 
\begin{Example}\label{crs}
  Let $(M,N,*)$ be a pair of based path connected topological spaces. Here
 $*\in N\subset M$. Then $\left (\pi_1(N,*),\pi_2(M,N,*), \d,\t \right
 )\doteq\Pi_2(M,N,*)$, where the boundary map and the action are the natural ones,  is a crossed module. See \cite{BV} or \cite{W3}- $IV.3$, for example. This is a result of J.H.C. Whitehead, see \cite{W1}.  
\end{Example}

A morphism $F=(\f,\p)$ between the crossed modules $\G$ and
$\G'=(G',E',\d',\t')$ is given by a pair of group morphisms $\f:G\to G'$ and
$\p:E \to E'$, making the diagram
\begin{equation*}
\begin{CD}
E @>\p>> E' \\
@V\d VV  @VV\d' V\\
 G @>>\f > G'
\end{CD}
\end{equation*}
commutative. 
In addition we must have:  $$\f(X)\t' \p (e)=\p(X \t e); \forall X \in G, \forall e \in E.$$
It is easy to prove that crossed modules and their morphisms form a category.
This category is a category with colimits, see \cite{BV}, $3.5$. 

\begin{Example}\label{base}
Let $(M,N,*)$ be a pair of based path connected spaces. If $*'\in N$ is another base point, and $\g$ is a path connecting $*'$ to $*$, then there exists a natural isomorphism $\Pi_2(M,N,*) \to \Pi_2(M,N,*')$, where the maps on fundamental groups are the usual ones constructed from the path $\g$.
\end{Example}

\begin{Example}\label{freeprod}
Let $\G=(G,E,\d,\t)$ and $\G'=(G',E',\d',\t')$ be crossed modules. The free product $\free(\G,\G')$ of $\G$ and $\G'$, also denoted by $\G \vee \G'$, is the pushout, in the category of crossed modules, of the diagram:
\begin{equation*}
\begin{CD}
(\{1\},\{1\},\t,\d) @>>> \G \\
@VVV  \\
\G'
\end{CD}.
\end{equation*}
Recall that the category of crossed modules is a category with colimits. This free product has the property that if $\cH$ is a crossed module then there exists a one to one correspondence between $\H\left (\G\vee \G', \cH\right)$ and $\H(\G,\cH) \times \H(\G',\cH)$.
We outline the construction of the free product of two crossed modules in example \ref{freecontruction}.
\end{Example}

\begin{Example}\label{simples}
Notice that if $E$ is an abelian group and the group $G$ acts on $E$ on the left by automorphisms, then $(G,E,\d=1_G,\t)$ is a crossed module.
Let $G$ be a group and $\kappa$ be a field. Let $E$ be the free $\kappa$ vector space on $G$, with its natural abelian group structure. Then $G$ has an obvious action $\t$ on $E$ on the left by linear transformations, thus $(G,E,\d,\t)$ is a crossed module. The explicit form of the action $\t$ is:
$$ X \t \left (\sum_{Y \in G} \lambda_Y .Y\right )=\sum_{Y \in G} \lambda_Y . XY \textrm{ where }  X \in G \textrm{ and }  \lambda_Y \in \kappa, \forall Y \in G.$$
This type of crossed modules is used in this article and \cite{FM} to prove the knottedness of some knotted surfaces.
\end{Example}

\begin{Example}\label{Relations}
Let $\G=(G,E,\d,\t)$ be a crossed module. Suppose that  the elements $a_1,...,a_n\in E$ are such that $\d(a_k)=1_G,k=1,...,n$. Let $F$ be the subgroup of $E$ generated by the elements of the form $X \t a_k$ where $X \in G,k=1,...,n$, thus $F$ is normal in $E$ by condition $2$ of the definition of crossed modules. Obviously both $\d$ and the action of $G$ on $E$ descend to $E/F$. Denote them by $\d'$ and $\t'$. It is easy to show that $(G,E/F,\d',\t')$ is a crossed module. We will go back to this important example in \ref{genrelations}.
\end{Example}

\subsubsection{Free Crossed Modules}\label{Refer10}
We now follow \cite{BV}, and mainly $3.1$.
Let $G$ be a group. There exists a subcategory of the category of crossed
modules for which the base group is always $G$ and the morphisms $F=(\f,\p)$ are such that $\f= \id_G$, see \cite{BV}, $2.2$. The objects of this category are  called crossed $G$-modules. Note that the subcategory of crossed $G$-modules  is not a full subcategory of the category of crossed modules.

\begin{Definition}
Let $K$ be a set and $\d_0:K \to G$ be a map.  A crossed module $\G=(G,E,\d,\t)$ is said to be the free crossed module (or free crossed $G$-module) on $K$ and $\d_0: K \to G$ if there exists an injective map $i:K \to E$ such that $\d \circ i=\d_0$, and the following universal property is satisfied:

\emph{
For any crossed $G$-module $\G'=(G,E',\d',\t')$ and any map $\p_0:K \to E'$ for which $\d' \circ \p_0=\d_0$, there exists a unique crossed $G$-module map $F=(\id_G,\p):\G \to \G'$ such that $\p \circ i=\p_0$.}

\end{Definition}
It is easy to prove that the free crossed modules just defined, if they exist, are unique up to isomorphism. See \cite{BV}. 

To prove their existence, let $E_0$ be the free group on the symbols $(X, m)$
where $X \in G$ and $m \in  K$. Therefore any element  $e \in E_0$ may be
expressed  as $e=(X_1,m_1)^{n_1} (X_2,m_2)^{n_2} ...(X_k,m_k)^{n_k}$, where
$X_i \in G$, $m_i \in K$ and $n_i \in \Z$, for $i\in \{1,...,k\}$. The group
$G$ acts on $E_0$ by automorphism as follows:  
\begin{multline*} Y \t \left ((X_1,m_1)^{n_1} (X_2,m_2)^{n_2} ...(X_k,m_k)^{n_k}\right )\\=(YX_1,m_1)^{n_1} (YX_2,m_2)^{n_2} ...(YX_k,m_k)^{n_k}.\end{multline*}
There also exists  a group morphism $\d:E_0 \to G$ such that: 
\begin{multline*}
\d((X_1,m_1)^{n_1} (X_2,m_2)^{n_2} ...(X_k,m_k)^{n_k})\\=\left (X_1\d_0(m_1)X_1^{-1}\right )^{n_1}\left (X_2\d_0(m_2)X_2^{-1}\right )^{n_2} ...\left (X_k\d_0(m_k)X_k^{-1}\right )^{n_k}.
\end{multline*}
The morphism $\d:E_0 \to G$ and the action $\t$ of $G$ on $E_0$ by automorphisms satisfy the condition $1.$ of the definition of crossed modules. However,  the pair $(G,E)$ provided with  $\d$ and $\t$ is in general not a crossed module, since the condition $2.$ of the definition of crossed modules may fail. In other words $(G,E_0,\d,\t)$ is, a priori, only a pre-crossed module, see \cite{BV}. To convert this construction into a crossed module, we define for any $e,f \in E_0$ the element:
$$[[e,f]]=efe^{-1}\left (\d(e) \t f^{-1}\right ),$$
called Peiffer commutators. It is easy to prove that  the subgroup $[[E_0,E_0]]$ of $E_0$ generated by all Peiffer commutators is a $G$-invariant normal subgroup of $E_0$ (this is true for any pre-crossed module). Indeed we have:
\begin{align*}
X\t [[e,f]] &= X\t \left (efe^{-1} \d(e) \t f^{-1}\right )\\
            &= \left ( X\t  e \right ) \left (X \t f\right) \left (X \t  e^{-1}\right )\left ( (X\d(e)) \t f^{-1}\right )\\
            &=\left ( X\t  e \right ) \left (X \t f\right) \left (X \t  e^{-1}\right )\left ( (X\d(e)X^{-1}) \t  (X \t f^{-1} )\right )\\ 
             &=\left ( X\t  e \right ) \left (X \t f\right) \left (X \t  e^{-1}\right )\left ( \d(X \t e)) \t  (X \t f^{-1} )\right )\\   
             &=[[X \t e,X \t f ]], \forall X \in G,\forall e,f \in E_0. 
\end{align*}
Moreover:
\begin{align*}
g[[e,f]]g^{-1}&=g\left (efe^{-1} \d(e)\t f^{-1} \right ) g^{-1}\\
              &=(ge) f (ge)^{-1} g    \d(e)\t (f^{-1}) g^{-1}\\
              &=(ge) f (ge)^{-1} \left (\d(ge) \t f^{-1}\right ) \left (\d(ge) \t f\right )g    \d(e)\t (f^{-1}) g^{-1}\\
              &=[[ge,f]] \left (\d(g) \t (\d(e) \t f) \right )g    \left (\d(e)\t f\right)^{-1} g^{-1}\\
              &=[[ge,f]][[g,\d(e) \t f]]^{-1}, \forall e,f,g \in E_0.
\end{align*}
In addition it is clear that $\d_{|[[E_0,E_0]]}=1_G$, since: $$\d\left (efe^{-1}\right )=\d(e)\d(f)\d(e)^{-1}=\d(\d(e) \t f),$$ the last equality follows from the fact that $(G,E_0,\d,\t)$ is a pre-crossed module.  Therefore $\d$ descends to a map $\d :E=E_0/[[E_0,E_0] \to G$.  For more details see \cite{BV}, section $3.3$. We have

\begin{Proposition}
In the notation just considered, $\G=(G,E,\t,\d)$ is a crossed module.
\end{Proposition} 
Notice that there exists a natural  map $i:K \to E$ for which $m \mapsto (1_G,m)$.  
\begin{Exercise}
Prove $i$ is injective.

\footnotesize{ {\bf  Hint} Consider a morphism $E_0 \to \Z$ for which $(X,m) \mapsto \lambda_m$, where $\lambda_m$ does not depend on $X$, and verify that it descends to $E=E_0/[[E_0,E_0]]$.} 
\end{Exercise}
It is trivial to conclude that with this map $\G$ is the free crossed module on $K$ and $\d :K \to G$. This is proved in \cite{BV}. We skip the proof since we are going to prove a more intricate result later. 
\begin{Remark}\label{reducerelations}
Actually it is possible to prove that $E\cong E_0/F$, where $F$ is the normal subgroup of $E_0$ generated by the elements of $E_0$ of the form:
$$\left (X,m \right) \left (Y,n \right)\left (X,m \right)^{- 1} \left (X \d_0(n) X^{-1}Y,n \right)^{-1}; X,Y \in G, m,n \in K.$$
See \cite{BV}, $3.3$ and \cite{BHu}.
\end{Remark}
We denote the free crossed module on $\d_0:K \to G$ by:
$${\cal F} \left (K \ra{\d_0} G\right ).$$

\subsubsection{Improving the Universal Property}\label{fund}
We freely use the notation of \ref{Refer10}.
Again let $G$ be a group and let $K$ be a set provided with a map $\d_0: K \to G$. Consider the free crossed module ${\cal F} (\d_0: K \to G)=(G,E,\t,\d)$  on $\d_0: K\to G$, thus $K$ is contained in $E$ under the injective map $i:K \to E$. This crossed module satisfies a stronger universal property than the one stated above. Let $G'$ be a group. Let also $\f: G \to G'$ be a morphism. Suppose we have a crossed module $\G'=(G',E',\d',\t')$. Consider a map $\p_0: K  \to E'$ such that $\d'\circ \p_0 =\f \circ \d_0$.

\begin{Theorem}\label{extension}
There exists a unique group morphism $\p:E \to E'$ extending $\p_0$ such that $(\f,\p)$ is a morphism of crossed modules.
\end{Theorem}  
\begin{Proof}
Recall the pre-crossed module $(G,E_0,\d,\t)$ constructed in \ref{Refer10}.
Consider the group morphism  $\p':E_0\to E'$ such that $\p'(X,m)=\f(X) \t' \p_0(m)$. Let us prove that $\f(X) \t' \p'(e)=\p'(X\t e); \forall e \in E_0, \forall X \in G$. We have:
\begin{align*}
\psi'\left ( X\t \left ( \prod _i (X_i,m_i)^{n_i} \right )\right )&=\psi'  \left ( \prod _i ( XX_i,m)^{n_i} \right)\\
&=  \prod _i (\f( XX_i)\t' \p'(m_i))^{n_i}\\&=
\f(X) \t' \left (\prod _i (\f(X_i)\t' \p'(m_i))^{n_i}\right )\\ &= 
\f(X) \t' \psi' \left ( \prod _i (X_i,m_i)^{n_i} \right ).\\ 
\end{align*} 
Let us now prove that $(\d' \circ \p')(e)=(\f \circ \d)(e), \forall e \in E_0$. We have

\begin{align*}
(\d' \circ \p')  \left ( \prod _i (X_i,m_i)^{n_i} \right )&=\d'  \left ( \prod _i (\f(X_i)\t' \p_0(m_i))^{n_i} \right )\\
&=   \prod _i \d' (\f(X_i)\t' \p_0(m_i))^{n_i} \\
&=   \prod _i  \left (\f(X_i)  (\d '\circ \p_0)(m_i) \f(X_i)^{-1}\right )^{n_i}\\
&=   \prod _i  (\f(X_i)  (\f \circ \d_0)(m_i) \f(X_i)^{-1})^{n_i}\\
&=  \f \left ( \prod _i  (X_i   \d_0(m_i) X_i^{-1})^{n_i}\right) \\
&=  \f \left ( \prod _i  (\d(X_i,  m_i) )^{n_i}\right) \\
&=(\f \circ \d )\left ( \prod _i (X_i,m)^{n_i} \right ). 
\end{align*}
With the equalities we prove that given $e,f\in E_0$ we have:
\begin{align*}
\p'(fef^{-1})&=\p'(f)\p'(e)\p'(f)^{-1}\\
&=(\d' \circ \p')(f) \t' \p(e)\\
&=(\f \circ \d)(f) \t' \p' (e)\\&=\psi'(\d(f) \t e).
\end{align*}
Thus $\p'$ descends to a group morphism $\p : E=E_0/[[E_0,E_0]] \to E'$. The rest of the argument is elementary.
\end{Proof}
\begin{Remark}
According to Ronald Brown, this result follows from the previous universal property in the category of crossed $G$-modules using the construction of induced crossed modules. See \cite{BV}, chapter $5$ and \cite{B2}.
\end{Remark}
\begin{Remark}
Notice that remark \ref{reducerelations} yields an alternative proof of theorem \ref{extension}.
\end{Remark}

\begin{Example}\label{free}
Let $X$ and $e$ be symbols. Let $G_X\cong \Z$ be the free group on $X$. The free crossed module $\G_{(X,e)}=(G_X,E_{(X,e)},\d,\t)$ on $X$ and $e$ is the free crossed module on $\d:\{e\} \to G_X$, where $\d(e)=X$. Thus in particular $E_{(X,e)}\cong\left <e\right >\cong\Z$ and $\d(e)=X$.
There exist inclusions $i_1:\{X\} \to G$ and $i_2:\{e\}\to E_{(X,e)}$, verifying $i_1(X) =\d \circ i_2(e)$.  This  crossed   module satisfies the following universal property: 

\emph{For any crossed module $\G=(G,E,\d,\t)$, and elements of $X_0\in G$ and $e_0 \in E$ such that $\d(e_0)=X_0$,  there exists a unique crossed module morphism $(\f,\p):\G_{(X,e)} \to \G$ such that $(\f \circ i_1)(X)=X_0$ and $(\p \circ i_2)(e)=e_0$.}

Therefore if $\G=(G,E,\d,\t)$ is a finite crossed module then: $$\#\H(\G_{(X,e)}, \G)=\#{E}.$$
\end{Example}

\begin{Example}\label{extension1}
Let $G$ be a group and $K$ be a set provided with a map $\d_0:K \to G$. Let $\G=(G,E,\d,\t)$ be the free crossed module ${\cal F }(\d_0:K \to G)$ on $\d_0:K \to G$, and let $i:K\to E$ be the inclusion map. Let $\G'=(G',E',\d',\t')$ be another crossed module. There exists a one-to-one correspondence between crossed module morphisms $(\f,\p):\G \to \G'$ and pairs $(\f,\p_0)$ where $\f:G \to G'$ is  a group morphism and $\p_0:K \to E'$ is a map verifying $\d'\circ \p_0=\f \circ \d_0$.   
\end{Example}

\subsubsection{Crossed Modules Presented By Generators and Relations}\label{genrelations}

Let $G$ be a group. Let also $K=\{m_1,...,m_k\}$ be a set, and let $\d_0:K \to G$ be a map. Consider the free crossed module ${\cal F}\left (\d_0: K \to G\right )=(G,E,\d,\t)$ on the map $\d_0:K \to G$. A relation on $(K,G)$ is, by definition, an element $r$ of $E$ with $\d(r)=1_G$. Since $E=E_0/[[E_0,E_0]]$, any  element of $E$ can be presented, though not uniquely in general, by a word of the form: 
$$r=(X_1,m_1 )^{\theta_1}...(X_l,m_l)^{\theta_l},$$
where $X_i \in G,\theta_i \in \Z$ and $m_i \in K$, $i=1,...,l$.
Let $R=\{r_1,...,r_n\}$ be a set of relations. The free crossed module on $\d_0:K \to G$ with  relations $\{r_1,...,r_n\}$, say ${\cal F}\left (\d_0:K \to G;r_1,...,r_n\right)$ is, by definition, the crossed module constructed from  ${\cal F}\left ( \d_0:K \to  G\right) $ and $\{r_1,...,r_n\}$ as in example \ref{Relations}. Note that we may very well take $R$  infinite.

 The isomorphism class of  ${\cal F}\left (\d_0:K \to G;r_1,...,r_n\right)$ is specified  by the universal property:

\emph{For any  crossed module $\G'=(G',E',\d',\t')$ and for any morphism $\f:G \to G'$, then given $e_1,...,e_k\in E'$ with $\d'(e_i)=(\f \circ \d_0)(m_i),i=1,...,k$, such that the induced map $(\f,\psi):{\cal F}\left (\d_0: K \to  G\right ) \to \G'$ verifies $\p(r_i)=1_{G'},i=1,...,n$, there exists a unique crossed module map $(\f,\psi '):{\cal F}\left (\d_0:K \to G;r_1,...,r_n\right)\to \G'$ such that $\p'(m_i)=e_i,i=1,...,k$.}

Therefore:

\begin{Theorem}\label{REFER}
Suppose that  $G$ is a free group, say on the set $L$, thus if $G'$ is a group then any map $\f_0:L\to G'$ determines, uniquely, a map $\f:G \to G'$. Let $K=\{m_1,...,m_k\}$ be a set provided with a map $\d_0: K \to G$. Let also $\{r_1,...,r_n\}$ be a set of relations. Let $\G'=(G',E',\d',\t')$ be a crossed module. There exists a one-to-one correspondence between crossed module maps:
 $${\cal F}\left (K \ra{\d_0} G;r_1,...,r_n\right)\to \G',$$
 and pairs of maps $\f_0:L\to G'$ and $\p_0:K \to E'$ verifying: 
\begin{enumerate}
\item  $(\f \circ \d_0)(m_i)=(\d' \circ \p_0)(m_i);i=1,...,k$,
\item $\psi(r_i)=1_{G'},i=1,...,n.$
\end{enumerate}
Here $(\f,\p)$ is the map ${\cal F}\left ( \d_0:K \to  G\right ) \to \G'$ determined by $\f_0$ and $\psi_0$. In particular if $r=(X_1, m_1 )^{\theta_1}...(X_l,m_l)^{\theta_l}$ then $\psi(r)=(\f(X_1)\t' \psi_0( m_1) )^{\theta_1}...(\f(X_l)  \t' \psi_0( m_l))^{\theta_l}\in G'$.
\end{Theorem}

\begin{Example}
Any crossed module $\G=(G,E,\d,\t)$ is the crossed module presented by $\d:E \to G$, with relations $(X,e)=(1,X \t e); X \in G, e \in E$ and $(1,e)(1,f)=(1,ef);e,f \in E$. 
 \end{Example}

\begin{Example}\label{freecontruction}
Let $\G=(G,E,\d,\t)$ and  $\G'=(G',E',\d',\t')$ be two crossed modules. The free product (defined in example \ref{freeprod}) $\G \vee \G'$ of $\G$ and $\G'$ is the crossed module over $G\vee G'$, on the map $E \times E' \to G \vee G'$, such that $(e,f) \mapsto \d(e)\d '(f)$,  considering the relations:
\begin{enumerate}
 \item $(X,e)=(1,X \t e); X \in G, e \in E$, 
 \item  $(X',e')=(1,X' \t e'); X' \in G', e' \in E$,
 \item  $(1,e)(1,f)=(1,ef); e,f \in E$,  
 \item  $(1,e')(1,f')=(1,e'f'); e',f' \in E'.$
\end{enumerate}
To prove this,  use the universal property of the free product of crossed modules. 
\end{Example}
Recall that the existence of the free product of crossed modules is a particular case of the theorem of existence of colimits in the category of crossed modules. This is proved, for example, in \cite{BV}.

\subsection{What is Wrong with $\pi_2(M,*)$ and a Solution }

As we mentioned before,  whenever  $E$ is an abelian group and the group $G$ acts on $E$ on the left by automorphisms, then $(G,E,\d=1_G, \t)$ is always a crossed module.
For any based path connected topological space $(M,*)$,  the group $\pi_2(M,*)$ is abelian and $\pi_1(M,*)$ acts on $\pi_2(M,*)$ by automorphisms. Therefore we have a crossed module $\pi_{1,2}(M,*)$ for any based topological space. See for example  \cite{L,M} for calculations of $\pi_{1,2}(M,*)$ where $M$ is the complement of a knotted surface in $S^4$.

A first idea concerning how to employ the notion of crossed modules to define
invariants of manifolds could be to consider the cardinal of morphisms from
the crossed module $\pi_{1,2}(M,*)$ into a crossed module $\G$. However, even
when $M$ is a compact manifold, it is not  certain that $\pi_2(M,*)$ is
finitely generated as a module over $\pi_1(M,*)$, thus there exists no
guarantee that this cardinal is finite, even when $\G$ is finite (c.f \cite{L}
problem $5$, and, less directly related, problems 6 and 13). Therefore $\pi_{1,2}(M)$  is not a very practical invariant since it not easy to distinguish two non-finitely  generated $\pi_1(M,*)$-modules.

Another solution is to consider the more tractable relative case. Let $(M,N,*)$ be a pair of path connected topological spaces with a base point $* \in N \subset M$. Therefore there exists a crossed module  
$\Pi_2(M,N,*)=\left (\pi_1(N,*),\pi_2(N,N,*),\d,\t\right)$, where all the structure maps are the natural ones (example \ref{crs}). As we will see, in the case when $N$  is the 1-skeleton $M^1$ of the connected CW-complex $M$ then $\Pi_2(M,M^1,*)$ is, in principle, easy to calculate, and finitely generated. In addition this crossed module contains all the information about the 2-type of $M$, therefore being stronger than $\pi_{1,2}(M)$ alone. See \ref{general}.

Choosing the apparently less charming relative case is also justified by the fact that $\Pi_2(M,M^1,*)$ does not depend on the cellular decomposition of $M$ up to free products with $\G_{(X,e)}=\Pi_2(D^2,S^1,*)$, see example \ref{free}. This permits us to obtain a homotopy invariant of finite connected CW-complexes for any finite crossed module, see \ref{mainset}.  
\subsubsection{A Theorem of Whitehead}
This is a basic theorem concerning the calculability of $\Pi_2(M,M^1,*)$ where $M^1$ is the 1-skeleton of the connected CW-complex $M$.
\begin{Example} 
\label{Whitehead} {\bf (Whitehead's Theorem)}
Let $M$ be a path-connected topological space with a base point $*$. Let $N$ be a topological space obtained from $M$ by attaching some 2-cells (or 2-handles) $s_1,...,s_n$. Choose a base point in the boundary of  any of the 2-cells, such that the attaching map $f_i:\d s_i \to M$ of any cell sends its base point to the base point $*$ of $M$, a fact we can suppose up to homotopy .  Each 2-cell $s_i$ therefore induces an element $\f(s_i)$ of $\pi_1(M,*)$, through its attaching map. Then, by a theorem of J.H.C. Whitehead in \cite{W}, the crossed module $\Pi_2(N,M,*)$ is the free $\pi_1(M,*)$-crossed module over the map $\f :\{s_1,...,s_n\} \to \pi_1(X,*)$. See \cite{BV}, $5.4$ and \cite{B,B2,GH}. 
\end{Example}
This is one of the most important results that we will use.

\begin{Example}\label{S1}
  Continuing the example above, a corollary of Whitehead's Theorem is  that $\Pi_2 \left (D^2,S^1,*\right)$, where $*\in S^1\cong \d D^2$,  is the free crossed module on the symbols $X$ and $e$ of example \ref{free}. Here $D^2=[0,1]^2$. Therefore if $\G=(G,E,\d,\t)$ is a finite crossed module then: $$\# \H\left (\Pi_2\left(D ^2,S^1,*\right ),\G\right)=\# E.$$ 
One can obviously prove this from the homotopy exact sequence of the pair $(D^2,S^1)$. 
\end{Example}

\subsubsection{The crossed module $\Pi_2(M,M^1,*)$ Is Finitely
 Generated}\label{finite}

As promised, we now prove this basic result:

\begin{Theorem}
Let $M$ be a finite connected CW-complex with a unique $0$-cell which we take as being its base point.  Let also $M^1$ be the $1$-skeleton of the cell decomposition of $M$. The  group
$\pi_2(M,M^1,*)$ is finitely generated as a $\pi_1(M^1,*)$-module. 
\end{Theorem}
\begin{Proof}
Let $M^2$ be the $2$-skeleton of $M$. By Whitehead's Theorem (example
\ref{Whitehead}), $\Pi_2(M^2,M^1,*)$ is the free  $\pi_1(M^1,*)$-crossed module
on the set of 2-cells of $M$ and their attaching maps. But since the set of
2-cells is finite, the construction of  free crossed $G$-modules in
\ref{fund} indicates that  $\pi_2(M^2,M^1,*)$ is finitely generated as a   $\pi_1(M^1,*)$-module. To finish the proof, we need to show that the inclusion of $(M^2,M^1,*)$ in
$(M,M^1,*)$ induces an epimorphism  $\pi_2(M^2,M^1,*) \to
\pi_2(M,M^1,*)$. This follows from the Cellular Approximation Theorem.  \end{Proof}

\begin{Corollary}
 Let $M$ be a finite connected CW-complex with a unique $0$-cell which we take as being its base point $*$. If $\G$ is  a finite crossed module, then the
number of morphisms from $\Pi_2(M,M^1,*)$ into $\G$ is finite.
\end{Corollary}
\begin{Proof}
This follows from the theorem above and the fact $\pi_1(M^1,*)$ is itself
finitely generated since $M$  is finite.
\end{Proof}

\subsubsection{2-Dimensional van Kampen Theorem}\label{Kampen}
It is important to note that the 2-dimensional van Kampen Theorem due to Ronald Brown and Philip Higgins is,  probably, the most useful tool to calculate the fundamental crossed modules of a based pair. For example, Whitehead's Theorem as well as theorem \ref{Attach3}, proved later, are  immediate corollaries of it.

The basis of this subsection is \cite{BV,BH}, explicitly theorem $2.3.1$ of \cite{BV} or theorem $C$ of \cite{BH}. See also \cite{B2}.
\begin{Theorem} {\bf (Brown and Higgins)}
 Let $(Y,X,*)$, where $* \in X \subset Y$, be a pair of based path connected topological spaces. Let $Y_1$ and $Y_2$ be path connected  subsets of  $Y$ whose union is $Y$. Suppose:
\begin{enumerate}
\item  the interiors of  $Y_1$  and $Y_2$ cover $Y$,
\item  the sets $X \cap Y_1$, $X\cap Y_2$, $Y_1 \cap Y_2$ and $X \cap  Y_1
  \cap Y_2$ are connected,  
\item  the inclusions of $X\cap Y_1$ in $Y_1$, of   $X\cap Y_2$ in $Y_2$, and of $X \cap (Y_1 \cap Y_2)$ in $Y_1 \cap Y_2$ induce epimorphisms of fundamental groups, 
\end{enumerate}
then the diagram
\begin{equation*}
\begin{CD}
\Pi_2(Y_1 \cap Y_2, X \cap Y_1 \cap Y_2,*) @>>> \Pi_2(Y_1 , X \cap Y_1,*)\\
@VVV                                              @VVV\\
 \Pi_2(Y_2 , X \cap Y_2,*) @>>> \Pi_2(Y_1 \cup Y_2=Y , X,*)
\end{CD},
\end{equation*}
with arrows induced by inclusions, is a pushout in the category of crossed modules.
\end{Theorem}

A corollary of the 2-dimensional van Kampen Theorem is the following fact which we will need later:

\begin{Example}\label{REFER1}
Let $M$ and $N$ be CW-complexes with  unique 0-cells, which we take as being their base points. We have:
$$\Pi_2((M,M^1,*) \vee (N,N^1,*)) \cong \Pi_2(M,M^1,*) \vee \Pi_2(N,N^1,*).$$ 
\end{Example}

\subsubsection{A General Theorem}\label{general}
Nothing here is strictly new. See for example \cite{W}.
Let $M$ be a CW-complex.  For any $k \in \N$, let $M^k$ denote the $k$-skeleton of $M$. For simplicity, suppose that $M$ has a unique $0$-cell, which we take as being our base point $*$. Suppose also that each $k$-cell of $M$ is attached through  a map $(S^{k-1},*) \to (M^{k-1},*)$, which we can suppose apart from homotopy, therefore avoiding complications related to the choice of base points. 
\begin{Remark}
Note that from  the Cellular Approximation Theorem we have $\Pi_2(M,M^1,*) \cong \Pi_2(M^3,M^1,*)$. 
\end{Remark}
Consider the group complex: 
$$... \ra{\d_4} \pi_3 (M^3,M^2,*) \ra{\d_3} \pi_2(M^2,M^1,*) \ra{\d_2} \pi_1(M^1,*)\ra{p} \pi_1(M,*),$$
which is a crossed complex of free type. See for example \cite{B2,B3}.
In particular, $\d_i \circ \d_{i+1}=1, \forall i$, and, moreover, for any $i \ge 2$, the group $\pi_1(M^1,*)$ acts on $\pi_i(M^i,M^{i-1},*)$ by automorphisms. If $i\ge 3$, then the action of $\pi_1(M^1,*)$ on $\pi_i(M^i,M^{i-1},*)$ is calculated from the action of $\pi_1(M^{i-1},*) \cong \pi_1(M)$  on  $\pi_i(M^i,M^{i-1},*)$ and $p:\pi_1(M^1,*) \to \pi_1(M,*)$, in the obvious way. Recall that all maps $\d_i$ are $\Z(\pi_1(M^1,*))$-module morphisms. 

Note that $\pi_2(M^2,M^1,*)$ is known by  Whitehead's Theorem (example \ref{Whitehead}), whereas $\pi_3 (M^3,M^2,*)$ is the free $\Z(\pi_1(M^2,*))$ module on the $3$-cells of $M$, a more known result (also due to J.H.C  Whitehead), see \cite{W3}, V.1.  The following is well known:
\begin{Lemma}
We have:
$$\pi_2(M,M^1,*)=\pi_2(M^2,M^1,*)/{\rm im}(\d_3).$$
\end{Lemma}
\begin{Proof}
This follows from the homotopy exact sequence of the triple $(M^3,M^2,M^1)$. Indeed, the following sequence is exact:
$$  \to \pi_3(M^3,M^2) \ra{\d_3} \pi_2 (M^2,M^1) \to \pi_2(M^3,M^1) \to \pi_2(M^3,M^2)\to ,$$ 
and, moreover $\pi_2(M^3,M^2) \cong \{1\}$.  Recall that we have $\pi_2(M,M^1,*)=\pi_2(M^3,M^1,*)$.
\end{Proof}

Let $c_1^3,...,c_{n_3}^3$ be the 3-cells of $M$. Each one of them defines an element of $\pi_3(M^3,M^2,*)$, thus $\d_3(c^3_i)\in \pi_2(M^2,M^1,*),k=1,...,n_3$. Let also $C^2=\{c_1^2,...,c_{n_2}^2\}$ be the set of 2-handles, where $c_k^2$ attaches along $\d_2(c_k^2) \in \pi_1(M^1,*),k=1,...,n_2$, therefore defining a map $C^2 \ra{\d_2} \pi_1(M^1,*)$. 

The following result appears in \cite{W}
\begin{Theorem}\label{Attach3}
The crossed module $\Pi_2(M,M^1,*)$ is the free crossed module on the set of 2-cells and their attaching maps, with one relation for each 3-cell of $M$. More precisely:
$$\Pi_2(M,M^1,*)\cong {\cal F} \left ( C^2 \ra{\d_2} \pi_1(M^1); \d_3(c^3_1),...,\d_3(c^3_{n_3})\right ).$$
Note that $\d_2 \circ \d_3=1$.
\end{Theorem}
An alternative proof can be done using the 2-dimensional van Kampen Theorem.

\begin{Proof}
By Whitehead's Theorem and the previous lemma, we only need to prove that ${\rm im} (\d_3)$ is the subgroup of $\pi_2(M^2,M^1,*)$ generated by the elements $X \t \d_3(c^3_i)$, where $X \in \pi_1(M^1,*)$ and $i=1,...,n_3$. This follows from the fact that $\pi_3(M^3,M^2,*)$ is the free $\pi_1(M^2,*)$ module on the set of 3-cells of $M$.
\end{Proof}
\begin{Remark} \label{ReferTrefoil}
If $\G=(G,E,\d,\t)$ is a crossed module, and $m \in \ker (\d)$ then $X \t m$ only depends on the projection  $p(X)$ of $X$ in $\mathrm{coker} (\d)$. Here $X \in G$. This is because $\ker (\d)$ commutes with all $E$.
\end{Remark}

This theorem proves that, in principle, if $M$ is a CW-complex and $*$ is a 0-cell, then the  crossed module $\Pi_2(M,M^1,*)$ can be calculated.  Despite the asymmetry introduced by choosing a particular 1-skeleton of $M$, this crossed module determines $\pi_2(M,*)$ and $\pi_1(M,*)$,  which fit in the exact sequence:
$$\pi_2(M)\cong \ker (\d) \ra{\iota} \pi_2(M,M^1) \ra{\d} \pi_2(M^1) \ra{p} \mathrm{coker}(\d)\cong \pi_1(M),$$
since $\pi_2(M^1)=\{0\}$. 
In fact $\Pi_2(M,M^1,*)$  also determines the  k-invariant $k(M,*) \in H^3(\pi_1(M,*),\pi_2(M.*))$. The  cohomology class $k(M,*)$ is determined from the classical correspondence between 3-dimensional group cohomology classes and crossed modules. See for example \cite{B2,ML,KB,JH}, thus  the crossed module $\Pi_2(M,M^1,*)$ determines the topological 2-type of $M$, see \cite{MLW}. We will go back to this issue in \ref{2Types}.

In fact, as mentioned before, something stronger is true. Namely we will prove that $\Pi_2(M,M^1,*)$ is well defined for any cellular space $M$ (i.e. independent on the choice of a cellular decomposition of $M$), up to free products with $\Pi_2(D^2,S^1,*) \cong \G_{(X,e)}$ (see example  \ref{free}), a result contained in  \ref{mainset}.

\subsection{The Main Set of Results}\label{mainset}
Let $(N,M)$  be a pair of CW-complexes such that the inclusion of $M$ in $N$ is a homotopy equivalence. Let $M^1$ and $N^1$ be, respectively, the 1-skeletons of $M$ and $N$. Suppose that $M$ has a unique 0-cell, which we take as being our base point $*$, so that both $M$ and $N$ are well pointed.

The group $\pi_1(M^1,*)$ is the free group on the set $\{d_1,...,d_m\}$ of 1-cells of $M$. There exist also $c_1,...,c_n \in \pi_1(N^1,*)$ such that $\pi_1(N^1,*)$ is the free group $F(d_1,...,d_m,c_1,...,c_n)$ on the set    $\{d_1,...,d_m,c_1,...,c_n\}$.

\begin{Theorem}
There exists a homotopy equivalence:
$$(N,N^1,*)\cong (M,M^1,*) \vee (D^2,S^1,*)^{\vee n}.$$                   
 \end{Theorem}

\begin{Proof} \footnote {This argument arose in a discussion with Gustavo Granja.}  
Since $M$ is a subcomplex of $N$, and $N$ is homotopic to $M$, it follows that $M$ is a strong deformation retract of $N$. By the Cellular Approximation Theorem, we can suppose, further, that there exists a retraction $r:N \to M$ sending $N^1$ to $ M^1$, and such that $r \cong \id_N$, relative to $M$.  In particular if $k \in \{1,...,n\}$ then we have $c_k r_*(c_k)^{-1}=1_{\pi_1 (N)}$, (though this relation does not  necessarily hold in $\pi_1(N^1,*)$).
Define a map $$f:(P,P^1,*) \doteq (M,M^1,*) \vee \bigvee_{k=1}^n (D^2_k,S^1_k,*)\to (N,N^1,*)$$ in the following way: First of all, send $(M,M^1,*)$ identically to its copy $(M,M^1,*)\subset (N,N^1,*)$. Then we can send each $S^1_k$ to the element $c_k r_*(c_k)^{-1}\in \pi_1(N^1,*)$. Since these elements are null homotopic in $(N,*)$, this map extends to the remaining 2-cells of $(P,P^1,*)$.

 Let us prove that $f:(P,P^1,*) \to (N,N^1,*)$ is a homotopy equivalence. It suffices to prove that $f:(P,*) \to (N,*)$ and $f_{|P^1}:(P^1,*) \to (N^1.*)$ are based homotopy equivalences.  See for example \cite{M2}, $6.5$. Notice that the  result proved there is also valid in the base case, as long as all the spaces considered are well pointed, which is the case.
 
 We first show that $f$ is an equivalence of homotopy $(P,*) \to (N,*)$.
Let $r':(P,*)\to (M,*)$ be the obvious retraction, thus $r' \cong \id_{(P,*)}$. We have $r \circ f\cong r \circ f \circ r' =r'\cong \id_{(P,*)}$. On the other hand $f \circ r=r \cong \id_{(N,*)}$.

We show now that $f$ is a  homotopy equivalence $(P^1,*) \to (N^1,*)$. It is enough to prove that the induced map $f_*: \pi_1(P^1,*) \to \pi_1(N^1,*)$ is an isomorphism. Note that  $\pi_1(P^1,*)$ is (similarly with $\pi_1(N^1,*)$) isomorphic  with the free group on the set $\{d_1,...,d_m,c_1,...,c_n\}$. The induced map on the fundamental groups has the form: 
$$f_*(d_k)=d_k, k=1,...,m, \textrm{ and }f_*(c_k)=c_kr_*(c_k)^{-1}, k=1,...,n.$$
 Notice that $r_*(c_k) \in F(d_1,...,d_m), k=1,...,n$. Consider the morphism $g$ of $F( d_1,...,d_m,c_1,...,c_n)$ on itself such that:
$$g(d_k)=d_k, k=1,...,m \textrm{ and } g(c_k)=c_kr_*(c_k), k=1,...,n.$$
 Therefore $(f_* \circ g )(d_k)=d_k, k=1,...,m$ and $(f_* \circ g)(c_k)=f_*(c_k r_*(c_k))=f_*(c_k) f_*(r_*(c_k))=c_kr_*(c_k)^{-1}r_*(c_k)=c_k,k=1,...,n$. On the other hand $(g \circ f_*)(d_k)=d_k, k=1,...,m$ and  $(g \circ f_*)(c_k)=g(c_k r_*(c_k)^{-1})=c_k r_*(c_k) g(r_*(c_k^{-1}))= c_k r_*(c_k) r_*(c_k^{-1})=c_k,k=1,...,n$. This  proves that $g^{-1}=f_*$, which finishes the proof.
\end{Proof}

\begin{Corollary}\label{referee}
Let $M$ and $N$ be homotopic CW-complexes with  unique $0$-cells, which we take as their base points $*$ and $*'$. There exists $m,n \in N$ such that:

$$\Pi_2\left((M,M^1,*)\vee (D^2,S^1,*)^{\vee n}\right) \cong \Pi_2\left ((N,N^1,*')\vee (D^2,S^1,*')^{\vee m}\right),$$
thus, in particular, from the 2-dimensional van Kampen Theorem\footnote{This result was suggested by a referee of a previous version of this article.} :
$$\Pi_2(M,M^1,*) \vee \Pi_2(D^2,S^1,*)^{\vee n} \cong \Pi_2(N,N^1,*')\vee \Pi_2(D^2,S^1,*')^{\vee m}.$$
\end{Corollary}

\begin{Proof}
There exists a CW-complex $P$ homotopic to $M$ and $N$ such that both $M$ and $N$ are embedded in $P$ as sub CW-complexes. The result follows from the previous theorem and example \ref{base}.
\end{Proof}
\begin{Theorem}\label{main0}
Let $M$ be a finite CW-complex with a unique $0$-cell, which we take as being our base point. Let $\G=(G,E,\d,\t)$ be a finite crossed module. The quantity:

$$I_\G(M)=\frac{\# \Hom\left (\Pi_2(M,M^1,*), \G\right)}{(\# E)^{b_1(M^1)}}$$
is finite, does not depend on the CW-decomposition of $M$ and is a homotopy invariant of $M$.
\end{Theorem}
Notice that we only need to suppose that the CW-complex $M^2$ is finite.

\begin{Proof}
The finiteness of $I_\G(M)$ was proved in \ref{finite}. let $N$ be a CW-complex with a unique $0$-cell $*'$ and homotopic to $M$. By the  previous corollary:
$$\Pi_2(M,M^1,*) \vee \Pi_2(D^2,S^1,*)^{\vee n} \cong \Pi_2(N,N^1,*')\vee \Pi_2(D^2,S^1,*')^{\vee m}.$$ 
In particular, from examples \ref{freeprod} and \ref{free} it follows that:
$$\Hom(\Pi_2(M,M^1,*),\G)\times E^{n} \cong  \Hom(\Pi_2(N,N^1,*),\G)\times E^{m}, $$
as sets. 
The result follows from the fact that we necessarily have $b_1(M^1)+n=b_1(N^1)+m$.
\end{Proof}

\begin{Remark}
As we have seen, it is in principle possible to calculate $\Pi_2(M,M^1,*)$  if $M$ is a CW-complex and $M^1$ is the 1-skeleton of it, where $\Pi_2(M,M^1,*)$ does not depend of the CW-decomposition of $M$ up to free product of $\Pi_2(D^2,S^1,*)$.  However, $\Pi_2(M,M^1,*)$ is defined by generator and relations,  and it  can be complicated to distinguish  two crossed modules presented in this way. Therefore our counting invariant $I_\G$ appears as a tool  to separate finitely presented crossed modules (up to free product of $\Pi_2(D^2,S^1,*)$).
\end{Remark}

\begin{Theorem}\label{main}
Let $M$ be a compact connected manifold with a handle decomposition with a unique $0$-handle. Choose a base point $*$ in the $0$-handle of $M$. Let $M^{(1)}$ be the handlebody made from the $0$ and $1$-handles of $M$, which we call the  1-handlebody of $M$. Let $\G=(G,E,\d,\t)$ be a finite crossed module. The quantity:
$$I_\G(M)=\frac{\# \Hom(\Pi_2(M,M^{(1)},*), \G)}{(\# E)^{b_1(M^{(1)})}}$$
is finite and it is a homotopy invariant of $M$.
\end{Theorem}

\begin{Proof}
It is well known that a handle decomposition of a manifold $M$  determines a topological space  $\hat{M}$ of the same homotopy type of $M$, with a CW-decomposition where each $i$-handle of the manifold $M$ generates an $i$-cell of the $CW$-complex $\hat{M}$, see \cite{RS} for example. Intuitively, $\hat{M}$ is obtained from $M$ by shrinking any $i$-handle to  an $i$-cell: its core. If $M^{(n)}$ is the $n$-handlebody of $M$, the above equivalence provides a homotopy equivalence $(M,M^{(n)},*) \cong (\hat{M},\hat{M^{(n)}},*)=(\hat{M},\hat{M}^n,*)$. See also \cite{H}, $6.4$. The results follows from the theorem above and the fact that a handle decomposition of a compact manifold is necessarily finite.
\end{Proof}

This theorem should be compared with the results of \cite{Y,P,P2,FM}. 

\subsubsection{Relation with Algebraic 2-Types}\label{2Types}
\begin{Definition}
A 2-type is a based path-connected topological space $(M,*)$ such that $\pi_k(M,*)=\{0\},\forall k>2$. If $M$ is a connected CW-complex with a base point $*$ which is a 0-cell, then $2T(M,*)$ is the based cellular space defined from $(M^3,*)$ by killing all the homotopy groups $\pi_k(M,*)$ of $M$ with $k>2$ in the usual way (see for example \cite{AH}, example 4.17). It is well known that $2T(M,*)$ does not depend on the CW-decomposition of  $M$ up to homotopy equivalence.   
\end{Definition}
If $M$ is a  connected CW-complex, then the CW-complex $2T(M,*)$ is called the 2-type of $M$, or, more commonly, the second Postnikov section of $M$.
\begin{Definition}
Let $(M,*)$ be a based path-connected topological space. Then the algebraic 2-type of $(M,*)$ is given by the triple $A2T(M,*)=\left (\pi_1(M,*),\pi_2(M,*),k(M,*)\right)$, where $k\in H^3(\pi_1(M,*),\pi_2(M,*))$ is the  $k$-invariant (or first Postnikov invariant). See \cite{ML,EML,MLW}. Recall that if $M$ and $N$ are well pointed connected CW-complexes then if $A2T(M,*)=A2T(M,*)$ then if follows that $2T(M,*)\cong 2T(N,*)$. 
\end{Definition}
Let $M$ be a connected  CW-complex. 
It is well known that $A2T(M,*)=A2T(2T(M,*))$. On the other hand $\Pi_2(M,M^1,*)=\Pi_2(M^3,M^1,*)=\Pi_2(2T(M,*),2T(M,*)^1),*)$. Therefore, it follows:
\begin{Theorem}
Let $\G$ be a finite crossed module. The homotopy invariant $I_\G(M)$, only depends on the  2-type $2T(M,*)$ of $M$. In particular $I_\G(M)$ only depends on the algebraic  2-type of $M$. 
\end{Theorem}

An interesting issue is therefore to determine how useful  $I_\G$ is to separate algebraic 2-types. 

\section{Application: 2-Knot Complements}
 For details on knotted surfaces, in particular movie presentation of them, we refer the reader to \cite{CS,CRS,CKS}. We work in the smooth category.

\subsection{A Type of Handle Decomposition of 2-Knot Complements} \label{handledec} 
Let $\S\subset S^4=\R^4 \cup\{\infty\}$ be a knotted surface, in other words a (locally flat) embedding of a 2-manifold $\S$ into $S^4$. We want to calculate $\Pi_2(M,M^{(1)},*)$, where $M = S^4\setminus \n(\S)$ and $M^{(1)}$ is the handlebody made from the 0 and 1-handles of a handle decomposition of $M$ (the 1-handlebody of $M$).  Here $\n(\S)$ is an (open) regular neighbourhood of $\S$.
As an application, we will also be able to  calculate $I_\G(S^4\setminus \n(\S))$, where $\G$ is a finite crossed module, see theorem \ref{main}.

\begin{figure}
\begin{center}
\includegraphics{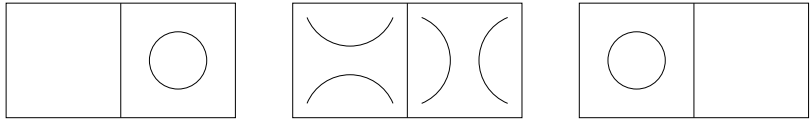}
\end{center}
\caption{Morse Modifications: Respectively a ``birth of a circle'',  a ``saddle point'' and  a ``death of a circle''.}
\label{crit1}
\end{figure}

\begin{figure}
\begin{center}
\includegraphics{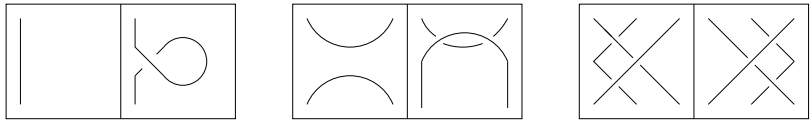}
\end{center}
\caption{Reidemeister moves I, II and III}.
\label{crit2}
\end{figure}

We construct, therefore, a natural handle decomposition of the complement of $\n(\S)$ in $S^4$:
 Suppose that  the projection on the last variable, the ``height'', is a Morse function on the knotted surface $\S$. In particular, for each non critical $t\in \R$, the set  $L_t =\S \cap \left (\R^3 \times \{t\}\right)$ is a link in $\R^3$ (a still of $\S$). Between critical values, the link $L_t$ will undergo an isotopy of $\R^3$. At critical points of degree $0$, $1$ or $2$, the link $L_t$ will suffer  Morse modifications, called, respectively, ``births of a circle'', ``saddle points'' and ``deaths of a circle''. See figure \ref{crit1}. A projection $p$ (which does not depend on $t$) onto some hyperplane of $\R^3$ will be a regular knot projection of $L_t$ for any non critical $t$, apart from a finite set. At these new non generic points the knot diagram will undergo a Reidemeister move change (see figure \ref{crit2}), and a planar isotopy will happen between them. The 1-parameter family of   projections $t \mapsto p (L_t)$ (or sometimes simply the 1-parameter family of links $t \mapsto L_t$), with the modifications at non-generic points,  will define what is  called a ``movie'' for the knotted surface $\S \subset S^4$. On the other hand  any movie defines, uniquely,  a knotted surface.

\begin{figure}
\begin{center}
\includegraphics{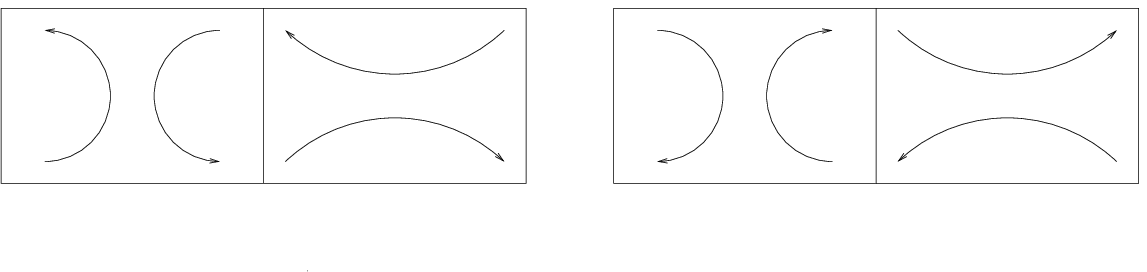}
\end{center}
\caption{Oriented saddle point moves.}.
\label{osaddle}
\end{figure}

Note that to the modifications displayed in figure \ref{crit2} (which may happen in any direction), one needs to add their mirror image. In the oriented case, one will need to consider all possible compatible orientations of all the diagrams. See figure \ref{osaddle} for the oriented version of saddle point moves.

Let $\S\subset  S^4$ be a knotted surface defined by a movie $t\mapsto L_t$. Choose  a projection $p$ of $\R^3$ onto a hyperplane of it so that $p(L_t)$ is a knot projection,  apart from a finite number of $t\in \R$. Therefore $\S$ is provided a Morse function and, away from critical points, $L_ t$ is a link in $S^3$. There exists a natural  handle decomposition of the complement of a regular neighbourhood of $\S$  where births/deaths of circles will induce $1/3$-handles of the decomposition, and saddle points induce $2$-handles. See \cite{GS}, section $6.2$, \cite{CKS}, $3.1.1$ or \cite{G}. This is very easy to visualise in dimension $3$. To calculate $I_\G(S^4\setminus \n (\S))$, however, we need an explicit description of this handle decomposition. We follow now \cite{CKS,G}, where  the missing bits of our description can be found.

\begin{figure} 
\begin{center}  
\includegraphics{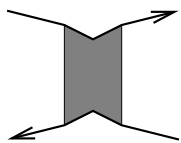} 
\end{center} 
\caption{A bit of a knot with bands.}
\label{kwb} 
\end{figure}

\begin{figure}
\begin{center} 
\includegraphics{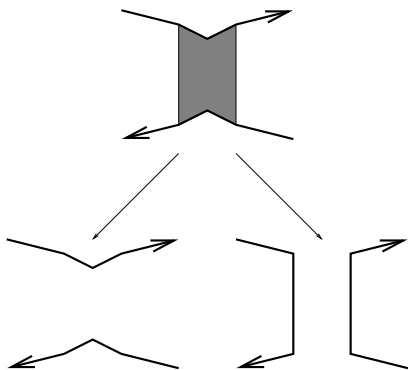} 
\end{center}
\caption{A piece of a knot with bands together with its pre- and post-knots (from
  left to right).}
\label{prepos}
\end{figure}

\begin{figure}
\begin{center}
\includegraphics{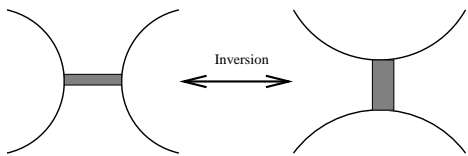}
\end{center}
\caption{A piece of a knot with bands and of the inversion of it.}
\label{inv}
\end{figure}

 What we are going to do is to present a movie for a  handle decomposition of
 the complement of a knotted surface $\S$, parallel to a (i.e. mounted upon a) given  movie of it. The 
 handle decompositions constructed are  understood more efficiently by considering knots with bands,
 possibly with spanners  \cite{SW} (this is the only point where we differ from \cite{CKS}). Recall that a knot with bands 
 is a knot together with some $I\times I$ bits, where $I=[0,1]$, intersecting
 the knot along $\d I \times I$. We assume that we have orientations on the knot as
 in  figure \ref{kwb}, or its mirror image.  A knot with bands $K$ determines two knots called the pre-knot and post-knot of $K$. See figure \ref{prepos}. The r\^ole of the  pre- and post-knots can be exchanged   by considering the inversion operation in knots with bands. See figure \ref{inv}. 

A spanner on a knot with bands $K$  is simply an unknot (or spanned component) of $K$ not linked to the pre-knot of $K$, together with a Seifert disk of it, which may intersect some  bands of $K$, transversally. These spanned unknotted components of $K$ are not considered to be part of the pre-knot of $K$.

\begin{figure}
\begin{center}
\includegraphics{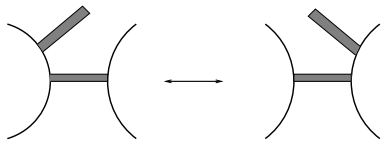}
\end{center}
\caption{Band Slide.}
\label{bandslide}
\end{figure}

\begin{figure}
\begin{center}
\includegraphics{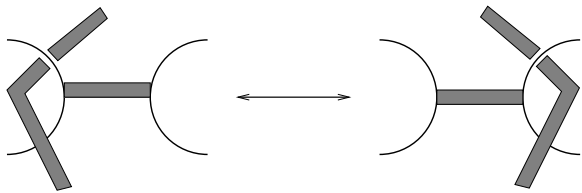}
\end{center}
\caption{Band Swim.}
\label{bandswim}
\end{figure}

\begin{figure} 
\begin{center} 
\includegraphics{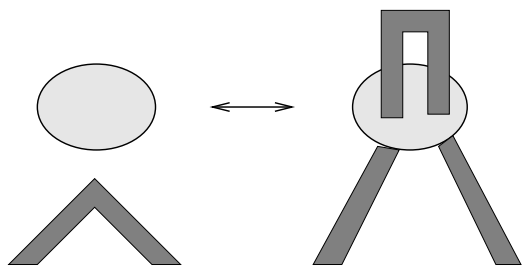} 
\end{center} 
\caption{A move on knots with bands and spanners. In both case they are good liftings of a knot in a movie to  a knot with bands and spanners at a maximal point.  This gives  us an example of how the attachment of $3$-handles is not uniquely defined.}
\label{ambiguous}
\end{figure}

As proved in \cite{SW}, section $4$, a movie $t \mapsto L_t$ of a knotted
surface $\S$ determines a movie $t \mapsto K_t$ of knots with bands and spanners,
such that the knot $L_t$ in each still is the pre-knot of $K_t$. These knot
with bands and spanners are defined up to the Band Swim Move and the Band Slide Move, depicted in figures
\ref{bandslide} and \ref{bandswim}, or, more precisely, their image under the inversion of knots with bands, and also the move in figure \ref{ambiguous}. See \cite{SW}. These extra components of the
movie of the knotted surface $\S$ only encode the previous stages of it, and should therefore, c.f. \cite{SW}, be considered as transparent. However this redundant information is very useful for our purposes.  Note that we can suppose that the (given) projection $p$ is a regular projection of $K_t\subset \R^3$, apart from a finite number of $t\in \R$.

We start by understanding the 1-handlebody $M^{(1)}$ of the handle decomposition of the
complement $M=S^4 \setminus \n(\S)$ of $\S$. If we have a birth of a circle, then the topology of the complement of $\S$ 
changes through the attachment of a 1-handle. Therefore, prior to the appearance
of saddle points, or deaths of a circle, the region of  the still (at some time $t$) of the movie
of the  handle decomposition of  $S^4 \setminus \n(\S)$ which corresponds to
its 1-handlebody (in other words the still of the $1$-handlebody at $t$) is the complement of a regular neighbourhood of the knot $L_t$ appearing in the still at $t$.

\begin{figure}
\begin{center}
\includegraphics{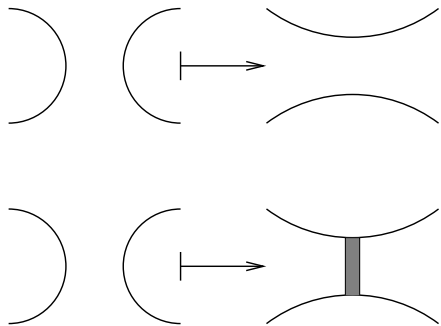}
\end{center}
\caption{Inserting bands on movies of knotted surfaces at saddle points.
\label{lift}}
\end{figure}

A saddle point determines a knot with bands $K$ such that the pre-knot of $K$
is the knot $L_{t_1}$ after the saddle point,  and the post-knot of it is the knot $L_{t_0}$
before the saddle point (notice the time inversion). See figure \ref{lift}. In all the forthcoming
stills of $t \mapsto K_t$, we need to consider the presence of this band, which co-evolves with the movie $t \mapsto L_t$ of $\S$ in the natural way (note that  between critical points $L_t$ undergoes isotopies of $\R^3$). See \cite{SW}, section $4$. When passing a saddle points, there are different possible configurations for a previously inserted band, but they are all related  by the Band Swim and Band Slide moves.

\begin{figure}
\centerline{\relabelbox 
\epsfysize 8cm
\epsfbox{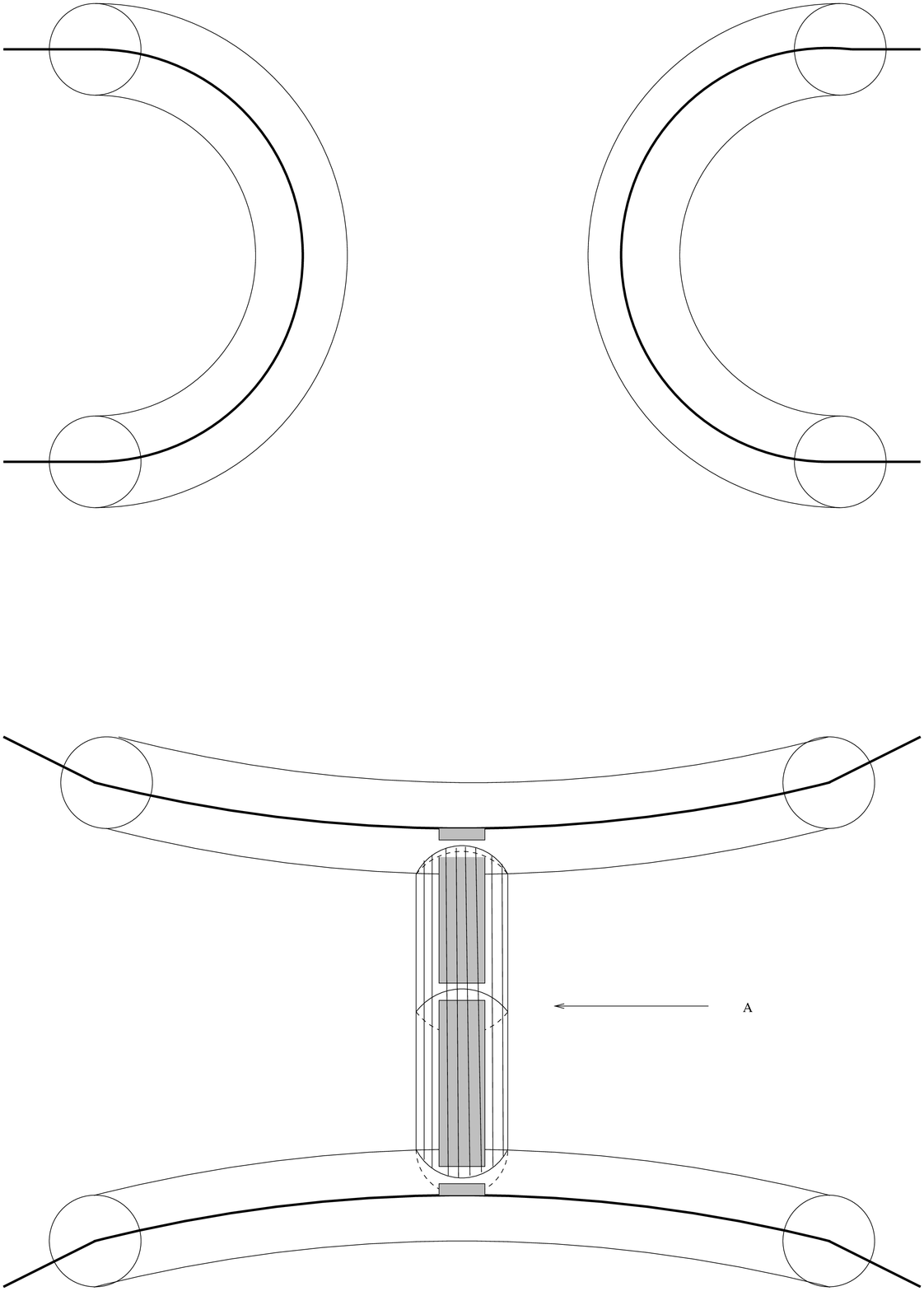} 
\relabel {A}{$t$-slice of 2-handle}
\endrelabelbox }
\caption{The Complement of a knotted surface after a saddle point, showing t-slice of attached 2-handle, constructed from the inserted band. The slice of its attaching region is the {\it lateral} border of the cylinder.}
\label{saddle}
\end{figure}

 Immediately after the saddle point, the configuration of the inserted 2-handle is as in figure \ref{saddle}. As time passes, the configuration can get more complicated, but it can always be read from the lifting of the movie $t \mapsto L_t$ to a movie $t \mapsto K_t$ of knots with bands and spanners, in the obvious way. Therefore, after the saddle point, the still of the $2$-handlebody of the handle decomposition of $S^4 \setminus \n(\S)$ will be precisely given by a regular
neighbourhood of the bands of $K_t$ minus a  regular neighbourhood of the pre-knot of $K_t$, which is $L_t$. The
bit of 1-handlebody left in this still will be given by the complement of a regular neighbourhood of the knot with bands. See  figure \ref{saddle}.

\begin{figure}
\begin{center}
\includegraphics{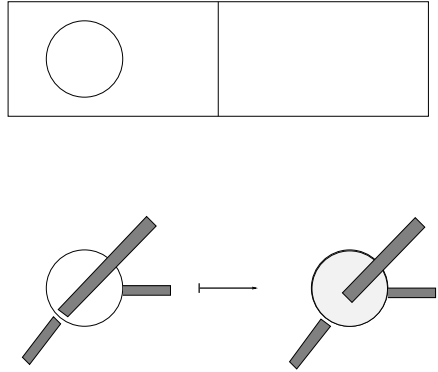}
\end{center}
\caption{Inserting spanners at deaths of a circle.
\label{maxspanner}}
\end{figure}

\begin{figure}
\begin{center}
\includegraphics{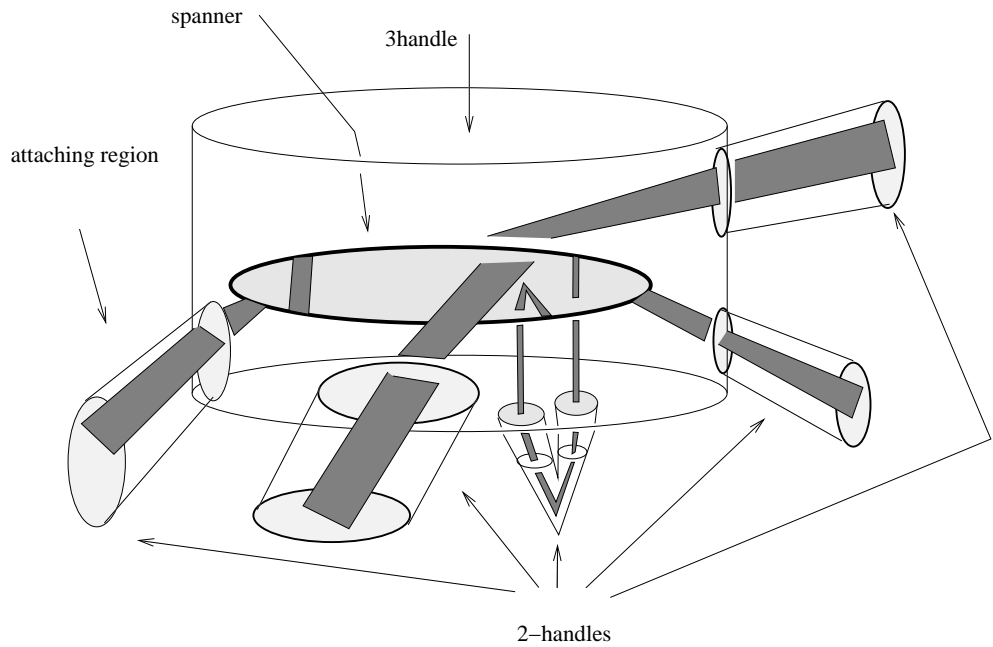}
\end{center}
\caption{A spanned component of a knot with bands defines a 3-handle of the
  complement of the knotted surface. Here the $t$-slice of the attaching region of
  the $3$-handle is all the boundary of the cylinder. Anything outside the solid regions belongs to the 1-handlebody $M^{(1)}$ of the handle decomposition. Note that some bands incident to the circle may not intersect the attaching region of the 3-handle.} \label{3handle}
\end{figure}

The (usually cumbersome) handles of index 3 attach along regions diffeomorphic with $S^2 \times I$. In the case of knot complements, $3$-handles correspond to deaths
of circles.  Let $t \mapsto L_t$ be a movie of $\S$. If we have a death of a circle at $t_0$, then in the forthcoming stills of $K_t$ (the lifting of $t \mapsto L_t$ to a movie of knots with bands and spanners) we will
have a knot with bands in which the circle still exists, but as a spanned
component, as in figure \ref{maxspanner}. See \cite{SW}, section $4$. The configuration of the attached 3-handle is as in figure \ref{3handle}. Therefore, for each forthcoming $t$, the still of the new
$3$-handle  will be given by a regular neighbourhood of the spanner, the still
of the 2-handlebody is given by a regular neighbourhood of the bands minus a
regular neighbourhood of all the spanners, and finally what is left of the
1-handlebody is given by the complement of a regular neighbourhood of all the
knot with bands, including all spanners.  

Finally we attach a 4-handle at $t_0$ when the link $L_t$  is empty for any $t\ge t_0$. 

 We have thus defined a handle decomposition of the complement of a knotted
 surface $\S$ if we are given  a movie $t \mapsto L_t$ of it. Notice that this handle
 decomposition  depends on the lifting of the movie $t \mapsto L_t$ to a movie $t \mapsto K_t$
 of  knots with  bands and spanners. Therefore the attachment of the 2-handles is  only defined up to
 handle slides, an ambiguity encoded by the Band Slide and Band Swim moves. The same is true for the $3$-handles, since the spanned components can have various shapes, and the same is true for their regular neighbourhoods; thus the attachment of 3-handles is also not uniquely defined. Therefore it needs to be treated case by case. See figure \ref{ambiguous}. These ambiguities have as a consequence that the usage of  Whitehead's theorem as well as the relations obtained from theorem \ref{Attach3} depend on the lifting of the  movie $t \mapsto L_t$ of a knotted surface $\S$  to a movie $t \mapsto K_t$ of knots with bands and spanners. These differences are quite obvious when we consider 3-handles.

\begin{figure} 
\begin{center} 
\includegraphics{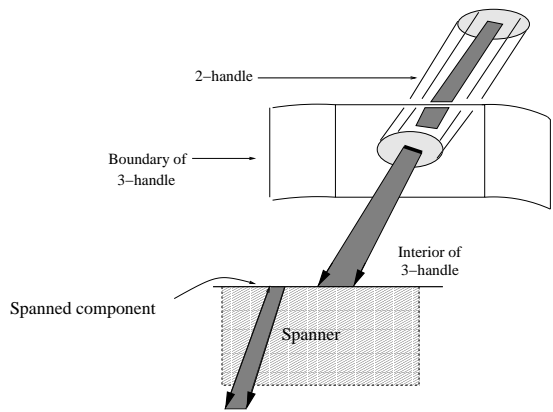} 
\end{center} 
\caption{Map of the attachment of a 3-handle determined by a spanned component in the vicinity of two bands. Note that only one of the bands intersects (a priori) the  attaching region of the 3-handle.}
\label{map}
\end{figure}

 In all the examples that we  will consider,  the spanned components will have an  obvious
 spanner, and the later obvious regular neighbourhoods. See figure \ref{map}.

\subsection{The Calculus}\label{calculus}
Let $\S$ be a knotted surface defined by a movie $t \mapsto L_t$, and let $t \mapsto K_t$ be a lifting of this movie to a movie of knots with bands and spanners. Choose a projection $p$ of $\R^3$ onto a hyperplane of it, so that $p(L_t)$ is a knot projection apart from a finite number $t$'s.
 We suppose that the surface $\S$ is oriented. Consider the handle decomposition of the 2-knot complement $M=S^4 \setminus \n(\S)$ just described,  and  let $M^{(1)},M^{(2)}$ and $M^{(3)}$ be the $1,2$  and $3$-handlebodies of it.  For any set $A$ in $S^4 =\R^4 \cup \{\infty\}$, let $A_t=A \cap (\R^3 \times \{t\})$ and $A_{\leq t}=A \cap (\R^3 \times \{(-\infty, t]\})$.

 In the representation of a knotted surface as a movie, we are free to take all births of circles at $t=-1$, all saddle points at $t=0$, and all deaths of circles at $t=1$. This kind of representation of knots are called hyperbolic splittings. For simplicity, we consider now movies of this type.

It is useful  to recall the calculation of the fundamental group of the complement of a knotted surface from a movie of it, explained, for example in  \cite{CKS}, $3.3.2$. Our construction can be seen as a perturbation of it. For the sake of simplicity, we forget about base points, and leave their trivial discussion to the reader. See example \ref{base}. For any $t\in \R$ we will choose a base point $*=*_t$ at the ``eye of the observer'' of the knot projection $p(L_t)$.

Let $t\in \R$ be such that the projection of $\S$ on the last variable (a Morse function) does not have
any critical point at $t$. We want to understand the crossed module
$\Pi_2(M_t,M^{(1)}_t)$.  Suppose $t \in (-1,0)$, thus, up to this point,  we only have $1$-handles in
the handle decomposition of the complement of $\S$. Therefore, in this case,
$\pi_1(M^{(1)}_t)=\pi_1(M^{(1)}_{\l})$ is simply the free group on the set of all
$1$-handles, and $\pi_2(M_t,M^{(1)}_{t})=\pi_2(M_\l,M^{(1)}_{\l})$ is the trivial group since $M_t =M^{(1)}_t$.

\begin{figure}
\begin{center}
\includegraphics{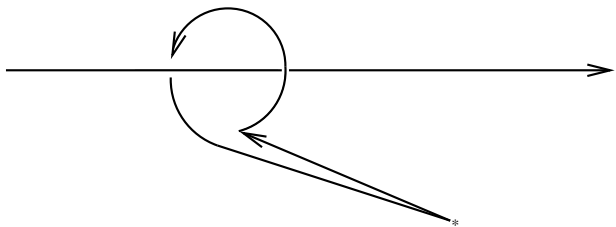}
\caption{Defining an element of the fundamental group of a knot complement. \label{fda}}
\end{center}
\end{figure}

\begin{figure}
\centerline{\relabelbox 
\epsfysize 3cm
\epsfbox{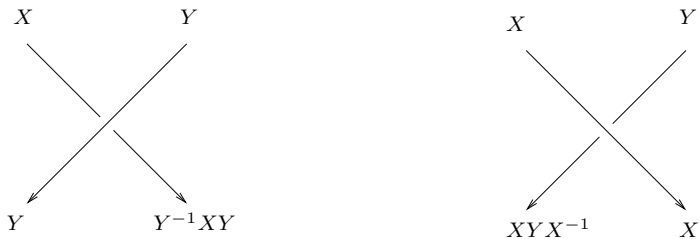} 
\relabel {X1}{$\s{X}$}
\relabel {X}{$\s{X}$}
\relabel {X3}{$\s{X}$}
\relabel {Y1}{$\s{Y}$}
\relabel {Y}{$\s{Y}$}
\relabel {Y2}{$\s{Y}$}
\relabel {Z}{$\s{Y^{-1}XY}$}
\relabel {W}{$\s{XYX^{-1}}$}
\endrelabelbox }
\caption{Wirtinger relations.}
\label{Colour}
 \end{figure}

 We
can define a representation of  $\pi_1(M^{(1)}_t)$ from  the Wirtinger
presentation of the fundamental group of the  complement of $L_t=K_t$. Therefore
each arc (upper crossing) of the projection $p(L_t)$ of $L_t$ generates an element of the fundamental group of
the complement of $L_t$, and each crossing generates a relation. See figures \ref{fda}
and \ref{Colour}. Notice that  we need to consider orientations on the knot diagram $p(L_t)$ 
so that these elements are well defined. Therefore, it is  at this point that
we need to introduce the (probably artificial) restriction that all the knotted surfaces that  we consider are oriented, so that we can have orientations on all $L_t, t \in \R$, `` evolving'' continuously with time.

\begin{figure}
\begin{center}
\includegraphics{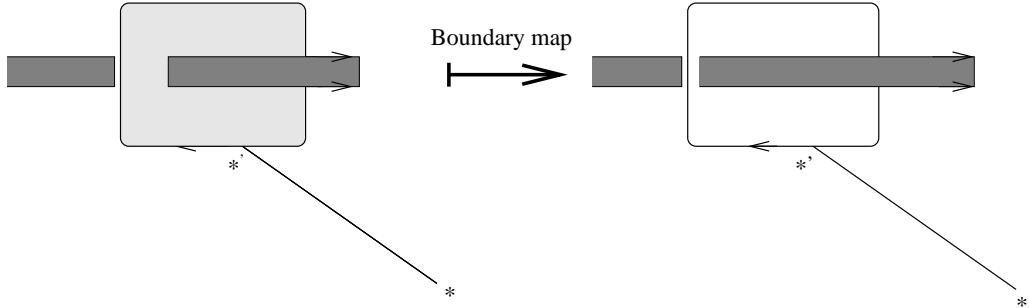}
\caption{A band (and a projection) defines an element of $\pi_2(M^{(2)}_t,M^{(1)}_t)$ as well as its boundary in $\pi_1(M^{(1)}_t)$.\label{fdb}}
\end{center} 
\end{figure}

When we pass the saddle points at $t=0$, we attach $2$-handles. Therefore, by Whitehead's Theorem (example \ref{Whitehead}), for any  $t \in
(0,1)$, the crossed module $\Pi_2(M_\l,M^{(1)}_{\l})$ is the free crossed module on the set of $2$-handles and their attaching maps, over the group $\pi_1(M^{(1)}_{<0}) \cong \pi_1(M^{(1)}_{\leq s}), s \in (0,1)$ (a free group). In particular,  the $2$-handles attached induce elements of $\pi_2(M_\l,M^{(1)}_{\l})$, as well as  elements of the group $\pi_1(M^{(1)}_{\l})$: their boundary. See figure \ref{fdb}.  

To specify the  element of $\pi_2(M_t,M^{(1)}_t)$ induced by each band (which will then induce an element of $\pi_2(M_\l,M^{(1)}_{\l})$), we will need to consider an orientation on each  band, as well as an  arc (upper crossing) of it in  projection $p(K_ t)$, a projection which is regular except for finite $t$.  The exact definition these elements of $\pi_2(M_t,M^{(1)}_t)$ is the following: Choose a regular neighbourhood of the band, which  is  thus diffeomorphic with $[0,1]^2 \times [0,1]$, where the last component is parallel to the core of the band. Choose a base point $*'$ on the surface of the regular neighbourhood, and a slice $[0,1]^2\times \{a\}$ of it incident to $*'$.  This slice thus defines an element of $\pi_2(M_t,M^{(1)}_t,*')$, as in figure \ref{fdb}. Suppose that $*'$ can be connected to the base point $*$ (the ``eye of the observer'' of $p(K_t)$) by a straight line which does not intersect the knot with bands, which  is thus totally contained in $M^{(1)}_t$. The natural morphism $\Pi_2(M_t,M^{(1)}_t,*')\to \Pi_2(M_t,M^{(1)}_t,*)$ defined by this curve determines the element of $\pi_2(M_t,M^{(1)}_t)$ specified by an arc of a band in $p(L_t)$.

\begin{figure}
\centerline{\relabelbox 
\epsfysize 5cm
\epsfbox{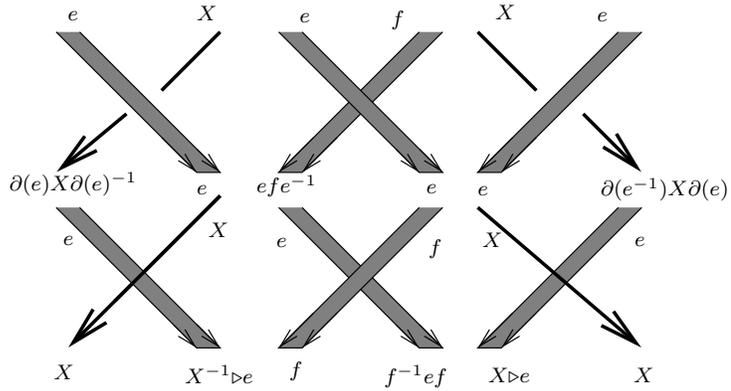} 
\relabel {A}{$\s{e}$}
\relabel {D}{$\s{e}$}
\relabel {B}{$\s{X}$}
\relabel {C}{$\s{\d(e)X\d(e)^{-1}}$}
\relabel {E}{$\s{e}$}
\relabel {F}{$\s{f}$}
\relabel {U}{$\s{efe^{-1}}$}
\relabel {V}{$\s{e}$}
\relabel {G}{$\s{X}$}
\relabel {H}{$\s{e}$}
\relabel {X}{$\s{e}$}
\relabel {Z}{$\s{\d(e^{-1})X \d(e) }$}
\relabel {I}{$\s{e}$}
\relabel {J}{$\s{X}$}
\relabel {K}{$\s{X}$}
\relabel {L}{$\s{X^{-1} \t e}$}
\relabel {M}{$\s{e}$}
\relabel {N}{$\s{f}$}
\relabel {O}{$\s{f}$}
\relabel {P}{$\s{f^{-1}ef}$}
\relabel {Q}{$\s{X}$}
\relabel {R}{$\s{e}$}
\relabel {S}{$\s{X\t e}$}
\relabel {T}{$\s{X}$}
\endrelabelbox }
\caption{Relations involving the elements of the relative second fundamental
  group induced by the bands of a knot with bands. Equation $2.$  of the
  definition of a crossed module is implicitly being used in the middle bits.}
\label{relations}
 \end{figure}

The group $\pi_1(M^{(1)}_t)$ is the fundamental group of the complement of a regular neighbourhood of the knot with bands and spanners $K_t$. This group can be calculated by the Wirtinger presentation for the fundamental group of graph complements, thus each arc of  $p(K_t)$ generates an element of  $\pi_1(M^{(1)}_t)$. As we have just seen, each arc a band in  $p(K_t)$ will  induce, further,  an element of $\pi_2(M^{(2)}_t, M^{(1)}_t)$. 
Further to Wirtinger relations, these group elements satisfy  the relations of figure \ref{relations}, easy to verify by a short calculation. 

\begin{figure} 
\centerline{\relabelbox 
\epsfysize 6cm
\epsfbox{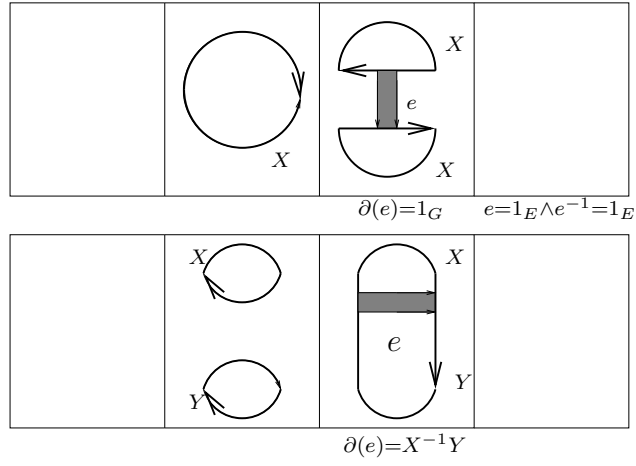}
\relabel {X}{$\s{X}$}
\relabel {Y}{$\s{X}$}
\relabel {Z}{$\s{X}$}
\relabel {rel1}{$\s{\d(e)=1_G}$}
\relabel {rel2}{$\s{e=1_E\wedge e^{-1}=1_E}$}
\relabel {e}{$\s{e}$}
\relabel {A}{$\s{X}$} 
\relabel {B}{$\s{Y}$}
\relabel {C}{$\s{X}$}
\relabel {D}{$\s{Y}$}
\relabel {f}{$e$}
\relabel {rel3}{$\s{\d(e)=X^{-1}Y}$}
\endrelabelbox }
\caption{\label{trivial} Generators and relations for  $\Pi_2(M_\l,M^{(1)}_\l)$, as  $t$ varies  in two examples. In both examples, the first  set of relations originates from the attachment of handles of index $2$ at $t=0$ (Whitehead's Theorem). In the first example, the second set of relations originates from the attachment of  $3$-handles at $t=1$. From figure \ref{map2}, for  the second example, no relation is implied from the attachment of a $3$-handle at $t=1$.} 
\end{figure}

Before explaining what happens after the remaining critical points, we
elucidate the calculation of the crossed module  $\Pi_2(M_\l,M^{(1)}_\l)=\Pi_2(M^{(2)}_\l,M^{(1)}_\l)$, where $t \in  (0,1)$. We
consider some simple examples.  We call the attention of the reader to figure
\ref{trivial}, the notation of which, is, we hope, self
explanatory. It follows  that we have, respectively in the first and second case:
\begin{align*}
\Pi_2(M^{(2)}_\l,M^{(1)}_\l)&={\cal F}\left (\{e\} \ra{ e \mapsto 1_G} F(X) \right),\\
\Pi_2(M^{(2)}_\l,M^{(1)}_\l)&={\cal F}\left (\{e\} \ra{ e \mapsto X^{-1}Y} F(X,Y)\right ),
\end{align*}
where $F(X)$ and $F(X,Y)$ are the free groups  on $\{X\}$ and on $\{X,Y\}$, respectively. Therefore if $\G=(G,E,\d,\t)$ is a finite crossed module then,  from example \ref{extension1}, we have:
\begin{align*}
  I_\G(M^{(2)}_\l)&=\frac {\#\{X\in G, e \in E |\d(e)=1_{G}\}}{\#E},\\
  I_\G(M^{(2)}_\l)&=\frac {\#\{X,Y \in G, e \in E |\d(e)=X^{-1}Y\}}{\#E^2},
\end{align*}
for the first and second cases, respectively. Here $I_\G$ is the invariant defined in theorem \ref{main}. Recall $t\in (0,1)$. 

\begin{figure} 
\centerline{\relabelbox 
\epsfysize 4cm
\epsfbox{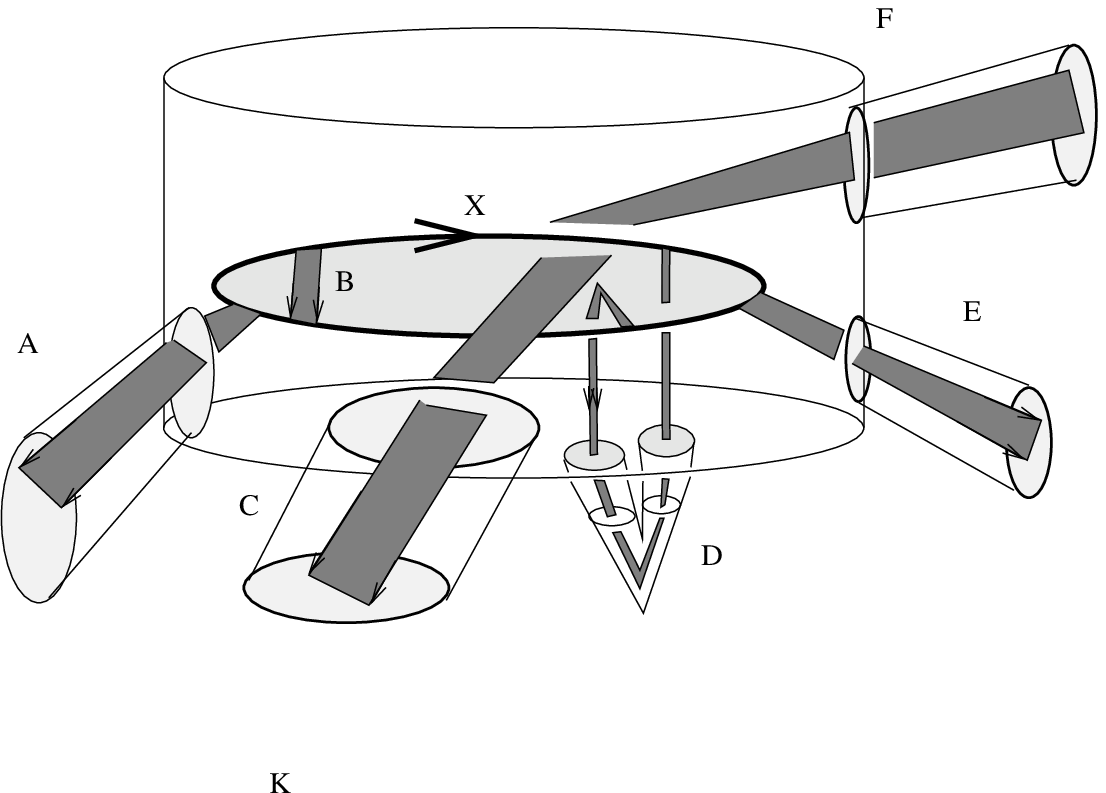}
\relabel {A}{$a$}
\relabel {B}{$b$}
\relabel {C}{$c$}
\relabel {D}{$d$}
\relabel {E}{$e$}
\relabel {F}{$X^{-1} \t c$}
\relabel {X}{$X$}
\relabel {K}{$\implies  a^{-1}c^{-1}d^{-1} d e^{-1}(X^{-1} \t c)=1$}
\endrelabelbox }
\caption{\label{3handlerel} Relations implied by the attachment of a 3-handle at a maximal point.} 
\end{figure}

Given the picture of the attaching region of a $3$-handle in figure
\ref{3handle},  it is easy to determine the relations (implied by theorem \ref{Attach3}) coming from the attachment 
of a $3$-handle at a maximal point: We simply need to determine the element $r$ of $\pi_2(M^{(2)},M^{(1)})$  induced by the attaching map of the 3-handle (see \ref{general}), and impose the relation $r=1$.  An example appears in figure \ref{3handlerel}. 
 Note that these relations depend on the spanner chosen for each disappearing circle, as well as on the regular neighbourhoods of it, since we need all this information to specify the  3-handle attached.

\begin{figure}
\begin{center}
\includegraphics{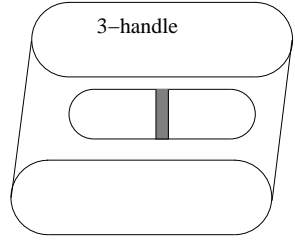}
\end{center}
\caption{Map of the attachment of the $3$-handle in the second example of figure \ref{trivial} and first example of figure \ref{trivial1}\label{map2}. Notice that the band shown  does not intersect the attaching region of the $3$-handle (the boundary of the cylinder). }
\end{figure}

 Two other examples of the use of the relations motivated by the attachment of 3-handles at maximal points are depicted in  figure \ref{trivial}.  In this figure, we choose  spanners for the disappearing circles as well as regular neighbourhoods of them  as in figure \ref{map}. For instance, in the second example of figure \ref{trivial}, the $3$-handle chosen has the shape shown in figure \ref{map2}, and, in particular, no relation appears from the attachment of  the $3$-handle at $t=1$. This is because, in  this case, the element of $\pi_2(M^{(2)},M^{(1)})$ induced by the attaching map of this 3-handle is      the identity of $\pi_2(M^{(2)},M^{(1)})$. See \ref{general}.

Therefore, we can conclude that in the first and second case of figure \ref{trivial}:
\begin{align*}
\Pi_2(M,M^{(1)})&={\cal F}\left (\{e\} \ra{ e \mapsto 1_G} F(X); e=1 \right),\\
\Pi_2(M,M^{(1)})&={\cal F}\left (\{e\} \ra{ e \mapsto X^{-1}Y} F(X,Y)\right ).\\
\end{align*}
Therefore if $\G=(G,E,\d,\t)$ is a finite crossed module then, we have, by theorem \ref{REFER}:
\begin{align*}
  I_\G(M)&=\frac {\#\{X \in G, e \in E |\d(e)=1_{G} , e=1\}}{\#E}\\
         &=\frac {\#\{X \in G, e \in E | e=1\}}{\#E} \\
         &=\frac{\#G}{\#E},
\end{align*}
and
\begin{align*}
  I_\G(M)&=\frac {\#\{X,Y \in G, e \in E |\d(e)=X^{-1}Y\}}{\#E^2}\\
         &=\frac{\#G}{\#E}.\\ 
\end{align*}
for the first and second cases, respectively.

\begin{figure} 
\centerline{\relabelbox 
\epsfysize 10cm
\epsfbox{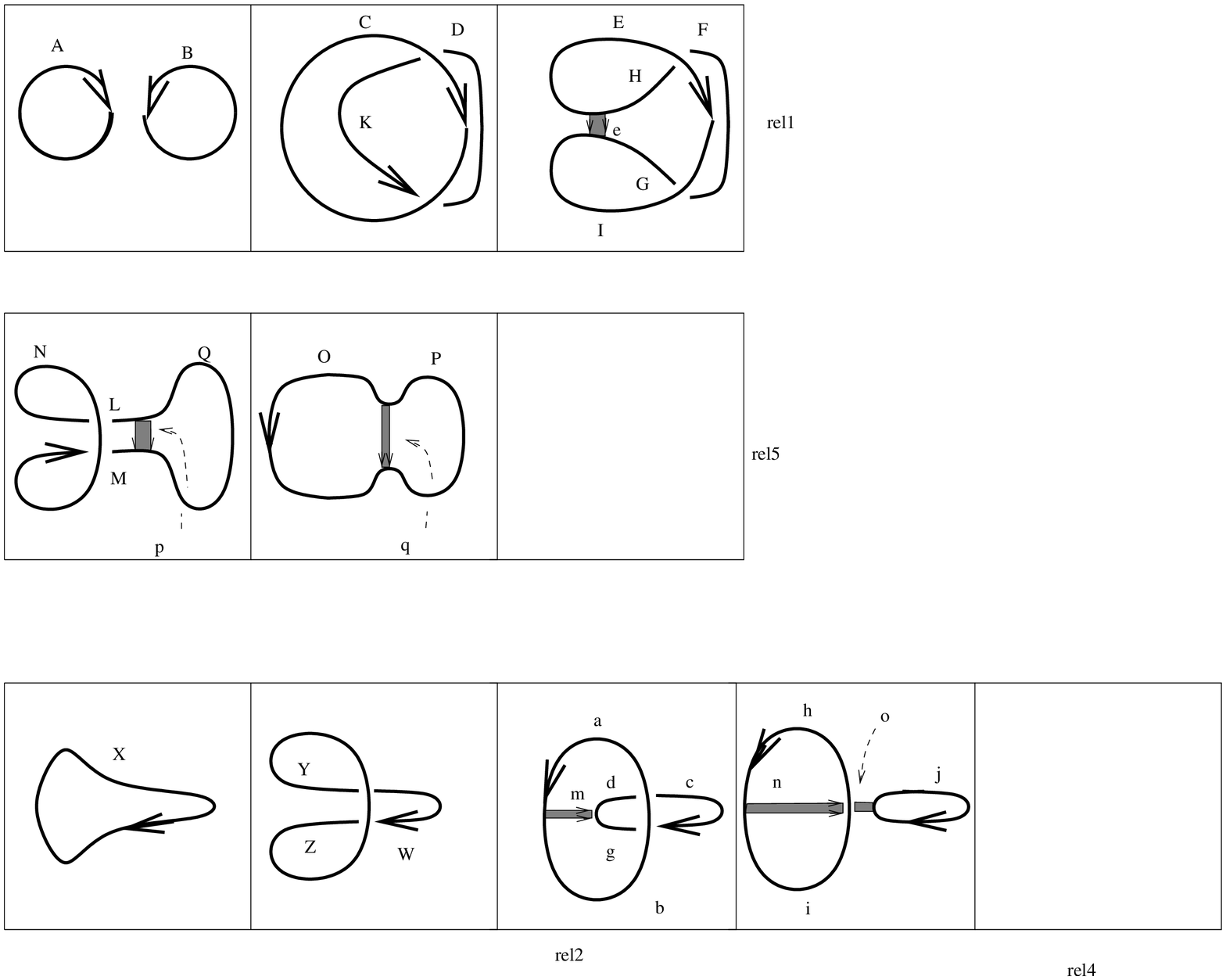}
\relabel {A}{$\s{X}$}
\relabel {B}{$\s{Y}$}
\relabel {C}{$\s{X}$}
\relabel {D}{$\s{Y}$}
\relabel {K}{$\s{Z}$} 
\relabel {E}{$\s{X}$}
\relabel {I}{$\s{X}$}
\relabel {F}{$\s{Y}$} 
\relabel {H}{$\s{Z}$}
\relabel {G}{$\s{Z}$} 
\relabel {e}{$\s{e}$}
\relabel {rel1}{\vbox{$\s{\d(e)=ZX^{-1}}$\\ $\s{Z=XYX^{-1}}$}}
\relabel {L}{$\s{X}$}
\relabel {M}{$\s{X}$}
\relabel {Q}{$\s{Y}$} 
\relabel {N}{$\s{X}$}
\relabel {O}{$\s{X}$}
\relabel {P}{$\s{Y}$}
\relabel {p}{$\s{X^{-1}\t e}$}
\relabel {q}{$\s{X^{-1} \t e }$}  
\relabel {X}{$\s{X}$}
\relabel {Y}{$\s{X}$}
\relabel {Z}{$\s{X}$}
\relabel {W}{$\s{X}$}
\relabel {a}{$\s{X}$}
\relabel {b}{$\s{X}$}
\relabel {c}{$\s{X}$}
\relabel {d}{$\s{X}$}
\relabel {g}{$\s{X}$}
\relabel {h}{$\s{X}$}
\relabel {i}{$\s{X}$}
\relabel {j}{$\s{X}$}
\relabel {m}{$\s{e}$}
\relabel {n}{$\s{e}$}
\relabel {o}{$\s{X \t e}$}
\relabel {rel2}{$\s{\d(e)=1_G}$}
\relabel {rel4}{$\substack{{X \t e=1_E}\\{X \t e^{-1}=1_E}}$}
\endrelabelbox }
\caption{\label{trivial1} Generators and relations for  $\Pi_2(M_\l,M_\l^{(1)})$ as $t$ varies: two  more examples. The first set of relations originates from Whitehead's Theorem. In the first example, from figure \ref{map2}, no relation originates from  the attachment of a $3$-handle at $t=1$. }
\end{figure}

 In figure \ref{trivial1}, we display the calculation of $\Pi_2(M,M^{(1)})$ in two other examples. We can conclude, that we have, respectively in the first and second cases: 
\begin{align*}
\Pi_2(M,M^{(1)})&={\cal F}\left (\{e\} \ra{ e \mapsto XYX^{-2}} F(X,Y) \right),\\
\Pi_2(M,M^{(1)})&={\cal F}\left (\{e\} \ra{ e \mapsto 1} F(X), X \t e=1\right ).
\end{align*}
Therefore, if $\G=(G,E,\d,\t)$ is a finite crossed module, then we have:
\begin{align*}
  I_\G (M)&=\frac {\#\{X,Y \in G, e \in E |\d(e)= XYX^{-2}\}}{\#E^2}\\      
          &=\frac{\#G }{\#E},
\end{align*}
and
\begin{align*}
 I_\G (M)&=\frac {\#\{X \in G, e \in E |\d(e)=1, X \t e=1\}}{\#E}\\
         &=\frac {\#\{X \in G, e \in E | e=1\}}{\#E}\\
         &=\frac{\#G }{\#E}.
\end{align*}    
for the first and second cases, respectively.

In figures \ref{trivial} and \ref{trivial1}, we depicted four examples of movies $t \mapsto L_t$ (or more precisely of their liftings $t \to K_ t$ to knots with bands) of the trivial embedding  of the sphere  $T$ with a  chosen orientation. These orientations are irrelevant for the final result, even though they are  needed so that we can consider the Wirtinger presentation of the fundamental group of the complement of $L_t$. Note that from all the examples it follows that $I_\G(T)=\#G/\#E$ if $\G$ is a finite crossed  module.

\begin{figure} 
\centerline{\relabelbox 
\epsfysize 6cm
\epsfbox{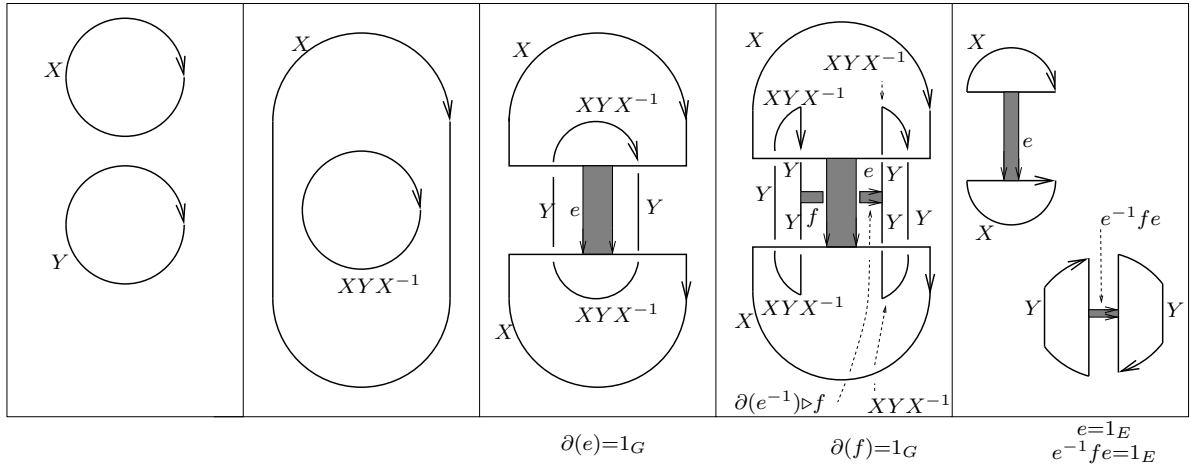}
\relabel{a}{$\s{X}$}
\relabel{b}{$\s{Y}$}
\relabel{c}{$\s{X}$}
\relabel{d}{$\s{XYX^{-1}}$}
\relabel{e}{$\s{X}$}
\relabel{f}{$\s{XYX^{-1}}$}
\relabel{g}{$\s{X}$}
\relabel{h}{$\s{XYX^{-1}}$}
\relabel{i}{$\s{Y}$}
\relabel{j}{$\s{Y}$}
\relabel{k}{$\s{e}$}
\relabel{rel1}{$\s{\d(e)=1_G}$}
\relabel{l}{$\s{X}$}
\relabel{m}{$\s{X}$}
\relabel{n}{$\s{XYX^{-1}}$}
\relabel{o}{$\s{XYX^{-1}}$}
\relabel{p}{$\s{Y}$}
\relabel{q}{$\s{Y}$}
\relabel{r}{$\s{Y}$}
\relabel{s}{$\s{Y}$}
\relabel{t}{$\s{XYX^{-1}}$}
\relabel{u}{$\s{XYX^{-1}}$}
\relabel{v}{$\s{Y}$}
\relabel{x}{$\s{Y}$}
\relabel{X}{$\s{e}$}
\relabel{Y}{$\s{f}$}
\relabel{Z}{$\s{\d(e^{-1})\t f }$}
\relabel{rel2}{$\s{\d(f)=1_G}$}
\relabel{A}{$\s{X}$}
\relabel{B}{$\s{e}$}
\relabel{C}{$\s{X}$}
\relabel{D}{$\s{Y}$}
\relabel{E}{$\s{e^{-1}f e}$}
\relabel{F}{$\s{Y}$}
\relabel{rel3}{$\substack{e=1_E\\{e^{-1}f e}=1_E}$}
\endrelabelbox }
\caption{\label{trivial3}Calculation of $I_\G(\S^4 \setminus \n(\S))$ where $\S$ is the trivial embedding of a disjoint union of two spheres. We skip the first and last stills, and consider the obvious spanning disk of the final unknotted components.}
\end{figure}

In figure \ref{trivial3} we display the calculation of $\Pi_2(M,M^{(1)},*)$, where $M=S^4 \setminus \n(\S)$, where   $\S$ is a trivial embedding of a disjoint union of two spheres.
This yields:
$$
\Pi_2(M,M^{(1)},*)={\cal F} \left ( \{e,f\}\ra{\begin{matrix}  e \mapsto 1 \\ f \mapsto 1 \end{matrix} } F(X,Y); e=1,e^{-1} f e=1\right ).
$$
As before we consider the obvious spanners for the unknots which disappear, as in figure \ref{map}. From this calculation it follows that $I_\G(S^4 \setminus \n(\S))=\frac{(\# G)^2}{(\# E)^2}$.

\begin{Remark}
Let $M$ be the complement of a knotted surface. Since $\pi_2(M^{(1)},*)=\{1\}$, the homotopy exact sequence of the based pair $(M,M^{(1)},*)$ permits to represent the second fundamental group of $M$ as: $$\pi_2(M,*)=\ker \left ( \d:\pi_2(M,M^{(1)},*) \to \pi_1(M^{(1)},*)\right ).$$
It is possible to  find a set of generators of $\pi_2(M,*)$. See \cite{BRS}.  
  Recall that $k(M)$ can  also be recovered from $\Pi_2(M,M^ {(1)},*)$, see \ref{general}.  Therefore our framework yields an algorithm to calculate $\pi_2(M,*)$ as a module over $\pi_1(M,*)$, as well as the first Postnikov invariant $k(M)$. 
 See \cite{L,M} for alternative approaches. 
\end{Remark}
\begin{Remark}
In \cite{PS}, it was asserted  the existence of pairs $\S,\S'\subset S^4$ of  knotted spheres such that $\pi_1(S^4\setminus \n(\S))=\pi_1(S^4 \setminus \n(\S'))$ and  $\pi_2(S^4\setminus \n(\S))=\pi_2(S^4 \setminus \n(\S'))$, as modules over $\pi_1(S^4\setminus \n(\S))=\pi_1(S^4 \setminus \n(\S'))$,  but with $k(S^4 \setminus \n(\S'))\neq k(S^4 \setminus \n(\S'))$. Since the invariant $I_\G(M)$ depends only on the algebraic $2$-type of $M$ (see \ref{2Types}), it would be interesting to determine whether $I_\G$ is strong, and practical,  enough to distinguish between pairs of embedded spheres with this property.  
\end{Remark}

\begin{Remark}
  Let $\G=(G,E,\d,\t)$  be a finite crossed module and let $M=S^4 \setminus \n(\S)$ be the complement of the knotted surface $\S$. To calculate $I_\G(M)$ we do not really need to determine $\Pi_2(M,M^{(1)},*)$ fully, which, since $\G$ verifies relations of its own, can be much more complicated than calculating $I_\G(M)$ only. By theorem \ref{REFER}, it is convenient to re-interpret the labellings introduced on the movies $t \mapsto K_t$ of  knots with bands  and spanners as the value on the generators of $\Pi_2(M,M^{(1)})$ of a morphism $\Pi_2(M,M^{(1)})\to \G$, where $\G=(G,E,\d,\t)$ is a finite crossed module, in which case we can add to the relations on figure \ref{relations} the relations satisfied in $\G$.    
\end{Remark}

\begin{Remark}\label{simple}
Continuing the last remark: Let $t \mapsto L_t$ be  a movie of a knotted surface $\S$ and let $t \mapsto K_t$ be a lifting of it to a movie of knots with bands and spanners. Let $p$ be a projection of $\R^3$ onto a hyperplane of it so that $p(K_t)$ is a regular projection, apart from a finite number of $t$. The  calculation of $I_\G(S^4 \setminus \n(\S))$  presented in this section simplify considerably when
$\G=(G,E,\d,\t)$ is a crossed module with $G$ abelian and $\d=1_G$, which
implies $E$ is abelian. In this case, for calculations purposes, at any time $t \in \R$, all components of the
knot with bands and spanners $K_t$  can pass through each other, except when one is a thin
component and the other is a band. In addition, the colourings of all the arcs
of a thin component of $p(K_t)$ are the same at all stages. See figures
\ref{Colour} and \ref{relations}, and adapt them to the case when $G$ is
abelian and $\d=1_G$.
\end{Remark}
\begin{Remark}
We defined our algorithm to calculate $\Pi_2(M,M^{(1)}.*)$ apriori only for hyperbolic splitting of knotted surfaces. However, it is easy to show that it remains valid for any movie of a knotted surface. We can prove this, for example,  by using the 2-dimensional van Kampen Theorem separately at each critical point, rather than using  Whitehead's Theorem followed by theorem \ref{Attach3}, as we did.
\end{Remark}
  
\subsection{Spun Hopf Link}

\begin{figure}
\begin{center}
\includegraphics{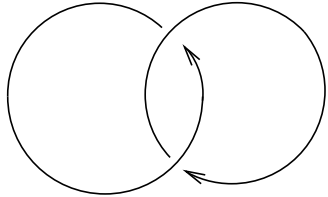}
\end{center}
\caption{Hopf Link.}
\label{hopflink}
\end{figure}

\begin{figure}
\begin{center}
\includegraphics{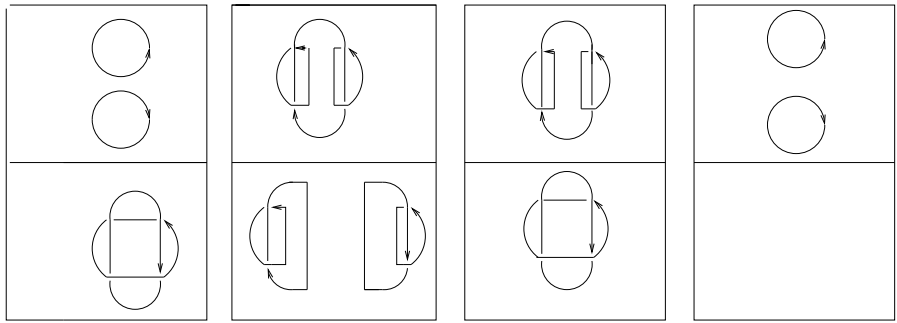}
\end{center}
\caption{Movie of the Spun Hopf Link. We skip the first stills and consider the
  movie of figure \ref{portion} between the first and second figures, and its
  inverse between the second last and last stills.}
\label{movie}
\end{figure}

\begin{figure}
\begin{center}
\includegraphics{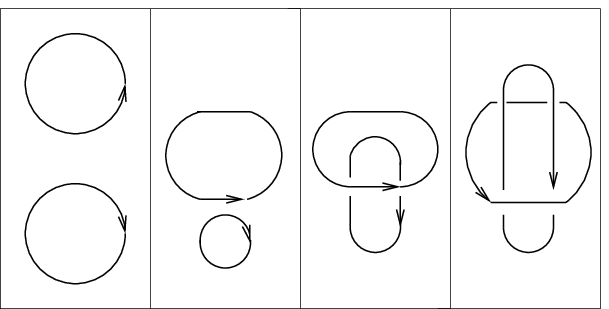}
\end{center}
\caption{A portion of the movie of the Spun Hopf Link.}
\label{portion}
\end{figure}

\begin{figure}
\begin{center}
\includegraphics{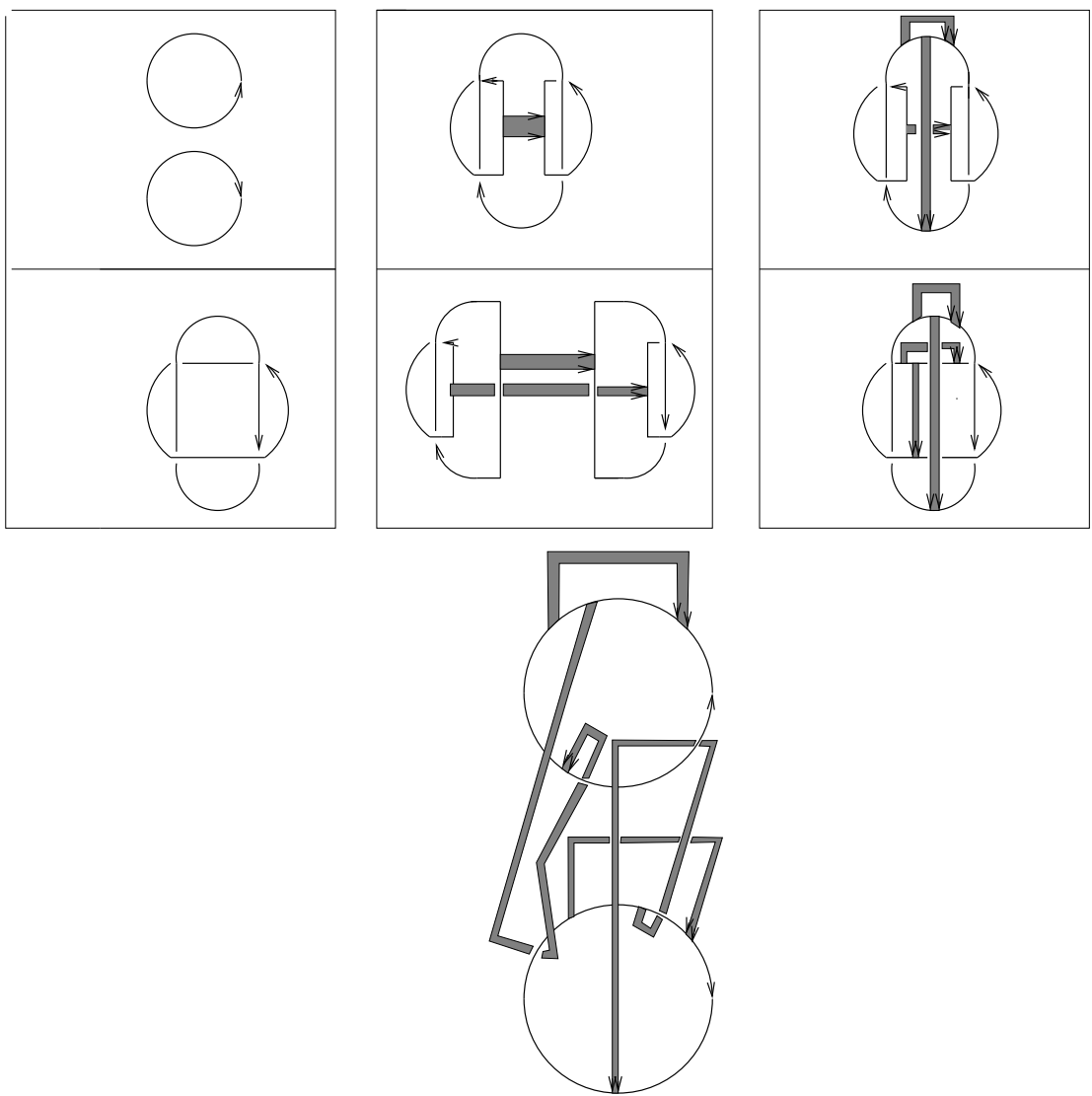}
\end{center}
\caption{Movie of the Spun Hopf Link, a version with bands.}
\label{movie2}
\end{figure}

We can now pass to a more complicated example. Consider the Spun Hopf Link
$\S$ obtained by spinning the Hopf Link depicted in figure \ref{hopflink}. A
movie for this knotted surface appears in figures \ref{movie} and
\ref{portion}. A lifting  of this movie of the Spun Hopf Link to a movie of
knots with bands and spanners  appears in figure \ref{movie2}. 

\begin{figure}
\begin{center}
\includegraphics{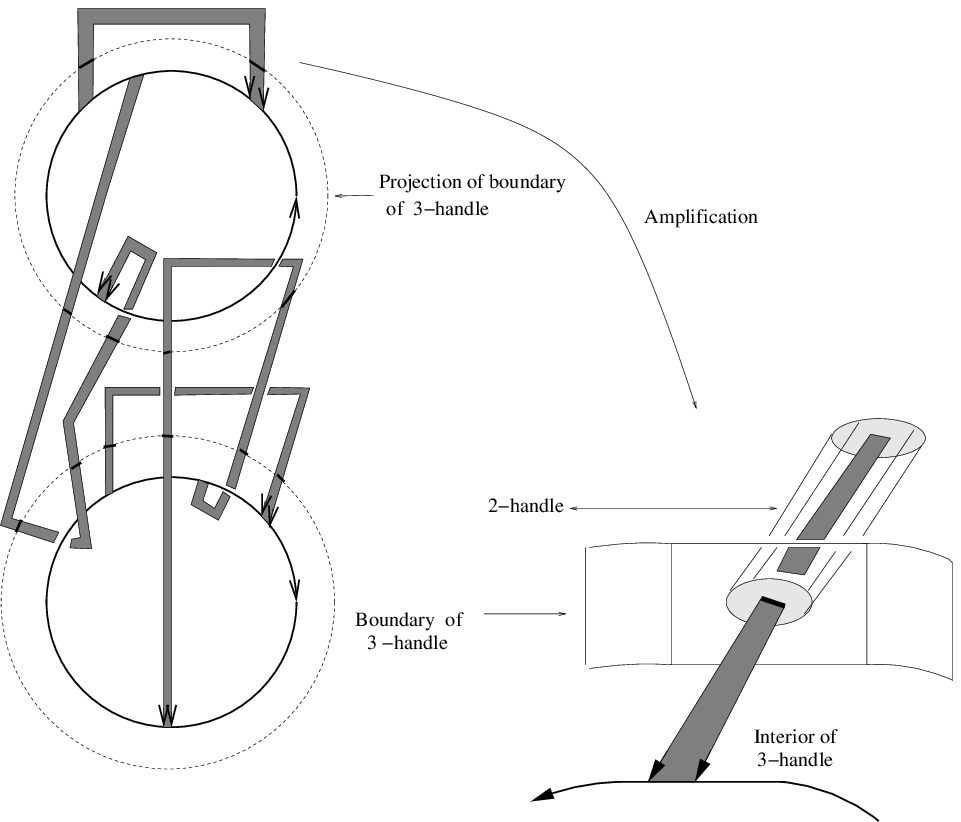}
\end{center}
\caption{Map of attaching region of $2/3$-handles in the last still figure \ref{movie2}.}
\label{upper}
\end{figure}

\begin{figure} 
\centerline{\relabelbox 
\epsfysize 6cm
\epsfbox{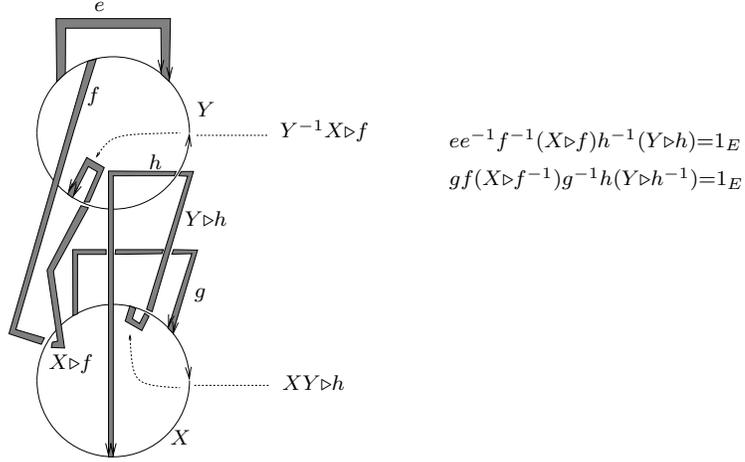}
\relabel {e}{$\s{e}$}
\relabel {f}{$\s{f}$}
\relabel {g}{$\s{Y^{-1}X \t f}$}
\relabel {h}{$\s{h}$}
\relabel {i}{$\s{Y \t h}$}
\relabel {j}{$\s{g}$}
\relabel {k}{$\s{X\t f}$} 
\relabel {X}{$\s{X}$}
\relabel {Y}{$\s{Y}$}
\relabel {m}{$\s{XY \t h }$} 
\relabel {rel}{\vbox{ $\s{ee^{-1}f^{-1}(X \t f) h^{-1} (Y \t h)} =1_E$\\ $\s{gf(X \t f^{-1}) g^{-1}h (Y \t h^{-1})=1_E}$}}
\endrelabelbox }
\caption{\label{Calc1} Calculation of $I_\G(S^4 \setminus \n(\S))$ where $\S$ is the Spun Hopf Link. Notice we suppose that $G$ is abelian and that $\d=1_G$. Otherwise several other labellings would need to appear.}
\end{figure}

Let $\G=(G,E,\d,\t)$ be a finite crossed module.
Using figures \ref{movie2} and \ref{upper}, as well as the relations of
figures \ref{Colour} and \ref{relations}, we can calculate $I_\G(S^4\setminus \n(\S))$. We will make the calculations only in the case in which $G$ is  abelian and $\d=1_G$, which implies that also $E$ is abelian. See remark \ref{simple}. We only show the last stage of the calculation, since the other steps are obvious. See figure \ref{Calc1}. 
From this we calculate:
\begin{equation}\label{value1}
I_\G(S^4\setminus \S)={\#\left \{X,Y \in G; f,h \in E: \substack{f^{-1} (X \t f) h^{-1} (Y \t h)=1_E\\   \textrm{and} \\f (X \t f^{-1})h (Y \t h^{-1})=1_E}  \right\}}
\end{equation}
We reiterate that this is only valid in the case in which $G$ is abelian and $\d=1_G$. Notice also our choice of attaching regions for $3$-handles which follows figure \ref{map}. Therefore, not every band incident to the circles will influence the calculation of the relations motivated by the attachment of $3$-handles.

On the other hand if $T$ is the trivial embedding of a disjoint union of two tori we have, for any crossed module $\G=(G,E,\d,\t)$:

$$I_\G(S^4 \setminus \n(T))=\left (\frac{(\#G) (\# \ker \d)^2}{\# E}\right )^2,$$
which specialises to: 
$$I_\G(S^4 \setminus \n(T))=(\#G)^2(\#E)^2, $$
in the particular case for which $G$ is abelian and $\d=1_G$.
\begin{Exercise}
Prove this.
\end{Exercise}
Let $E$ be an abelian group and $G$ be a group with an action $\t$  on $E$ by automorphisms. Then if we put $\d=1_G$ we can define a crossed module $\G=(G,E,\d,\t)$. Let $G=\Z_2$ and $E$ the free $\Z_2$ vector space on the set $\Z_2$ with the obvious action of $\Z_2$ by linear transformations. See example \ref{simples}. Then it follows that $\G=(G,E,\t,\d=0)$ is a crossed module. Note we have switched to additive notation, more adapted to this example. Let $\Z_2=\{[0],[1]\}$. If we choose $f,h=[1]$, $X=[1]$ and $Y=[0]$ we have $-f+ (X \t f)- h+  (Y \t h)\neq 0_E$. In particular it follows that the Spun Hopf Link is knotted.

 Notice that we would have not been able to prove this (known) fact if we had used
 only the fundamental group of the complement of the  Spun Hopf Link, and considering the number of
 morphisms from  it into an abelian group, since the first homology group of
 the complement only depends on the intrinsic topology of the knotted
 surface, and not on the embedding. Therefore we argued that homotopy 2-types are good invariants of
 complements of knotted surfaces. For a more detailed discussion see \cite{FM}.

\subsection{Explicit Calculation for the Spun Trefoil}

\begin{figure}
\centerline{\relabelbox 
\epsfysize 5cm
\epsfbox{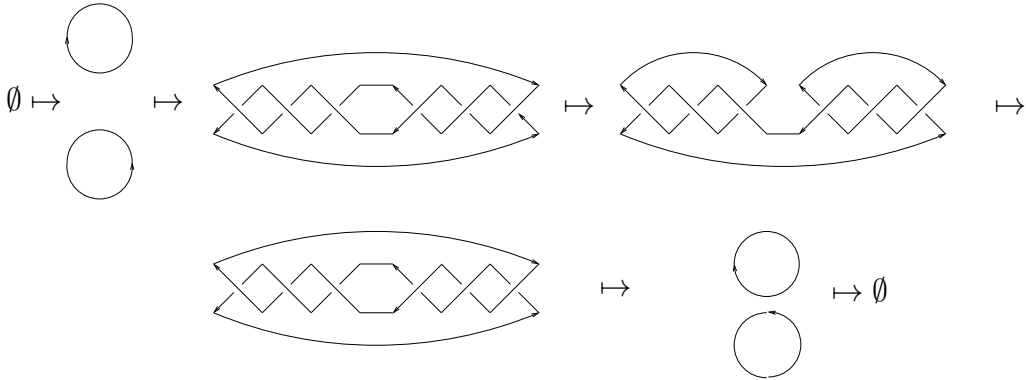}
\relabel{t1}{$\emptyset\mapsto$}
\relabel{t2}{$\mapsto$}
\relabel{t3}{$\mapsto$}
\relabel{t4}{$\mapsto$}
\relabel{t5}{$\mapsto$}
\relabel{t6}{$\mapsto\emptyset$}
\endrelabelbox}
\caption{ \label{ST} A movie of the Spun Trefoil.}

\end{figure}

A movie of the Spun Trefoil appears $\S$  in figure \ref{ST}. We assign the obvious
movies to the second and the second last arrows. Observe that the choice of
these two movies is irrelevant for the final result, anyway.

\begin{figure}
\centerline{\relabelbox 
\epsfysize 17 cm
\epsfbox{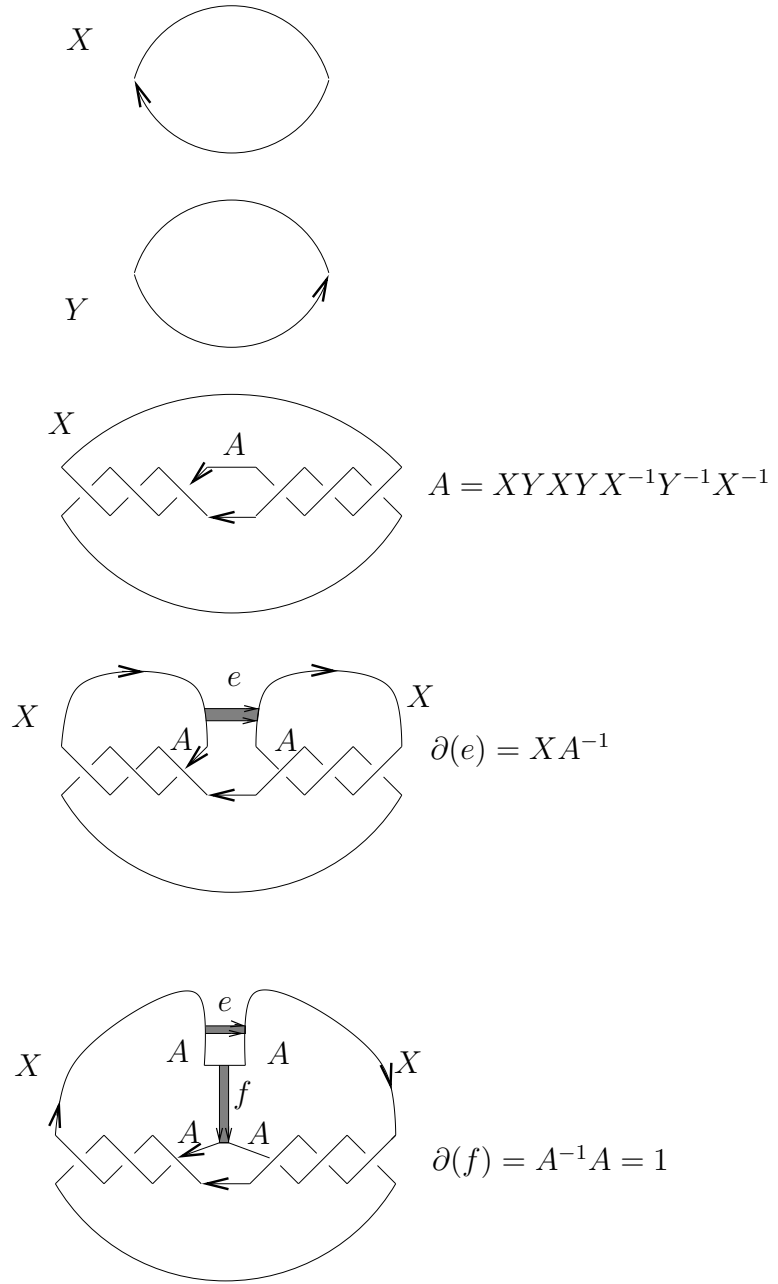}
\relabel{X0}{$X$}
\relabel{X1}{$X$}
\relabel{X2}{$X$}
\relabel{X3}{$X$}
\relabel{X4}{$X$}
\relabel{X5}{$X$}
\relabel{Y0}{$Y$}
\relabel{A0}{$A$}
\relabel{A1}{$A$}
\relabel{A2}{$A$}
\relabel{A3}{$A$}
\relabel{A4}{$A$}
\relabel{A5}{$A$}
\relabel{A6}{$A$}
\relabel{e0}{$e$}
\relabel{e1}{$e$}
\relabel{f0}{$f$}
\relabel{R1}{$A=XYXYX^{-1}Y^{-1}X^{-1} $}
\relabel{R2}{$\d(e)=X A^{-1}$}
\relabel{R3}{$\d(f)=A^{-1}A=1$}
\endrelabelbox} 
\caption{\label{STbands} Calculation of $\Pi_2(M,M^{(1)})$, where $M$ is the complement of the Spun Trefoil,  first part. }
\end{figure}

\begin{figure}
\centerline{\relabelbox 
\epsfysize 17cm
\epsfbox{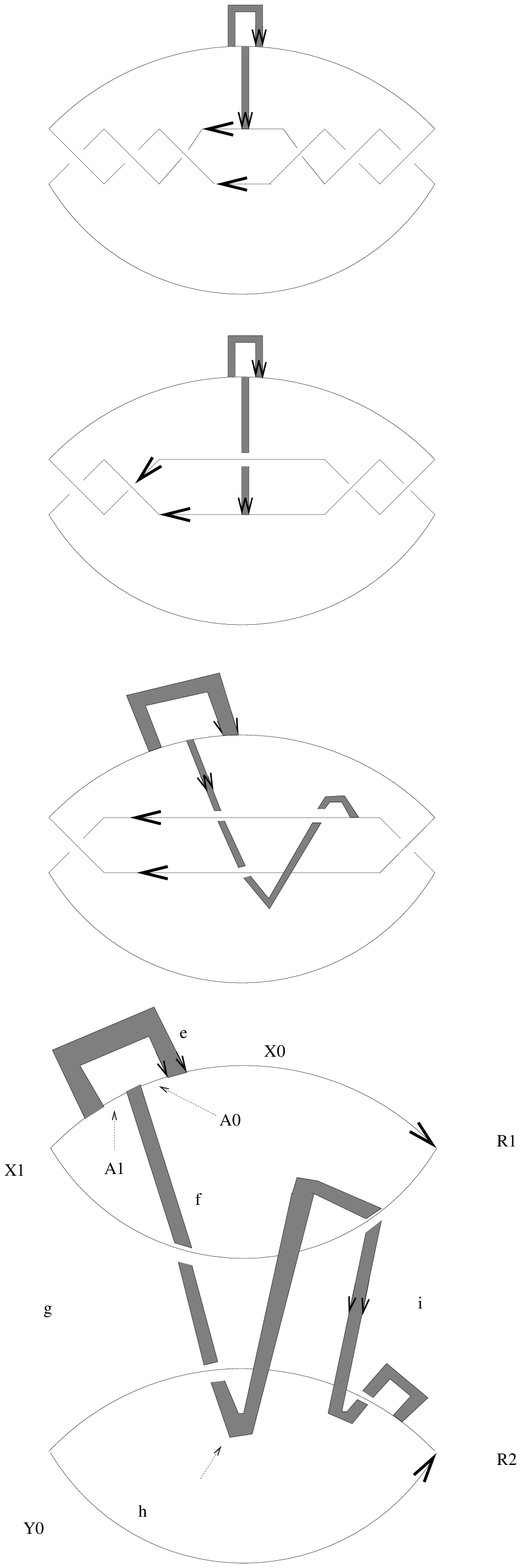}
\relabel{X0}{$X$}
\relabel{X1}{$X$}
\relabel{A0}{$A$}
\relabel{A1}{$A$}
\relabel{e}{$e$}
\relabel{f}{$f$}
\relabel{g}{$X^{-1} \t f$}
\relabel{h}{$Y^{-1}X^{-1} \t f$}
\relabel{i}{$X^{-1} Y^{-1}X^{-1} \t f$}
\relabel{Y0}{$Y$}
\relabel{R1}{$\s{(X^{-1}\t f^{-1})( Y^{-1} X^{-1} \t f)( X^{-1}Y^{-1} X^{-1} \t f^{-1})=1}$}
\relabel{R2}{$\s{(X^{-1}\t f^{-1})( Y^{-1} X^{-1} \t f)( X^{-1}Y^{-1} X^{-1} \t f^{-1})=1}$}
\endrelabelbox}
\caption{\label{ST2bands}Calculation of $\Pi_2(M,M^{(1)})$, where $M$ is the complement of the Spun Trefoil,  second part. }

\end{figure}

Let $M=S^4 \setminus \n(\S)$ be the complement of the Spun Trefoil $\S$, with the handle decomposition defined from the shown movie of $\S$  as in \ref{handledec}. 
We display in figures \ref{STbands} and \ref{ST2bands} the calculation of $\Pi_2(M,M^{(1)})$, following the method proposed  in \ref{calculus}. Similarly with the calculations for the Spun Hopf Link, we use the configuration of the attached 3-handles at the two maximal points as in figure \ref{map}.

This permits us to conclude that $\Pi_2(M,M^{(1)})$ is the free crossed module over $F(X,Y)$ on the map $\{e,f\}\to F(X,Y)$ such that $e\mapsto  XA^{-1}$, where $A=XYXYX^{-1}Y^{-1}X^{-1}$, and $f\mapsto 1$, considering the relation $(X^{-1}\t f^{-1})( Y^{-1} X^{-1} \t f)( X^{-1}Y^{-1} X^{-1} \t f^{-1})=1$.  Here $F(X,Y)$ is the free group on $X$ and $Y$.
From this it follows that, if $\G=(G,E,\d,\t)$ is a finite crossed module, with $G$ and $E$ abelian, and $\d=1$, then:
$$I_\G(M)=\frac{\#\left \{(X,f) \in G \times E:( X\t f)( X^2\t f^{-1 })( X^{3} \t f)\right \}}{\#E}.$$
This agrees also with the calculation in \cite{FM}. As proved there, this information suffices to prove that the Spun Trefoil is knotted. Indeed, if $G=\Z_3$ and $E$ is the free $\Z_2$ vector space on $\Z_3$ with the natural action $\t$  of $G$ and $\d=1$, then $I_\G(S^4 \setminus \n(T))\neq I_\G(S^4 \setminus \n(\S))$, where $T$ is the trivial embedding of $S^2$ in $S^4$ and $\G=(G,E,\d,\t)$.

\subsubsection{The Second Fundamental Group of the Spun Trefoil Complement}

The explicit construction of the free crossed modules and of a free crossed module with relations will be needed now.

\begin{figure}
\centerline{\relabelbox 
\epsfysize 3 cm
\epsfbox{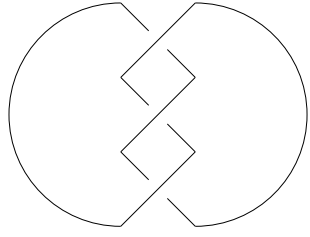}
\endrelabelbox}
\caption{Trefoil Knot.}
\label{trefoil}
\end{figure}

Let $K\subset S^3=\R^{3} \cup \{\infty\}$ be a knot. Suppose that  the projection on the last variable is a Morse function in $K$. Then, similarly with the 4-dimensional case, we have a handle decomposition of the complement of $K$ where local minimums induce 1-handles of the complement and local maximums induce 2-handles. We have at the end to attach an extra 3-handle, which will cancel out  one of the 2-handles previously  attached. Therefore, the complement of the Trefoil Knot, shown in figure \ref{trefoil}, admits a handle decomposition with one 0-handle, two 1-handles and one 2-handle. See \cite{GS}, exercise $6.2.2$. Note that 2-dimensional CW-complexes with a unique $2$-cell are classified by their fundamental group, up to homotopy equivalence. see \cite{J}.

Consider the CW-complex $N$ with one 0-cell $\{*\}$, two 1-cells $X$ and $Y$ and a 2-cell $e$ attaching along $XA^{-1}$. Here $A=XYXYX^{-1}Y^{-1}X^{-1}$. Then $N$ is homotopic to the Trefoil Knot complement. In particular $\pi_2(N,*)=\{0\}$, by the well known theorem (due to  Papakyriakopoulos)  asserting that 3-dimensional knot complements are aspherical.  See \cite{Pa}. On the other hand  we can represent $\pi_2(N,*)$ as $ \ker \{\d: \pi_2(N^2,N^1,*) \to \pi_1(N^1,*)\}$. Note that $\Pi_2(N^2,N^1,*)$ is the free crossed module on the map $e \mapsto XA^{-1}\in F(X,Y)$. 

Let again $M$ be the complement of the Spun Trefoil. Let $F$ be the subgroup of $\pi_2(M,M^{(1)},*)$ generated by the elements $B \t f$, where $B \in F(X,Y)$. Since $\d(f)=1$ this group commutes with all of $\pi_2(M,M^{(1)},*)$ and, in particular, any element of $\pi_2(M,M^{(1)},*)$  expresses, uniquely, as $g=ab$, where $a \in F$ and $b\in \pi_2(N,N^1,*)$. Thus, since $\d(F)=1$, we may have $\d(g)=1$ if, and only if, $\d(b)=1$ or  $b=1$\footnote{This is where the asphericity of 3-dimensional knot complements is being used.}. In particular, $\pi_2(M,*)=\ker \{\d: \pi_2(M,M^{(1)},*) \to  \pi_1(M^{(1)},*)\}=F$. Let $G=F(X,Y)/(XA^{-1}=1)$, which is isomorphic with the fundamental group of the complement of the Trefoil Knot.  By remark \ref{ReferTrefoil}, it follows that $F$ is the free abelian module over $G$ with one generator $f$ and the relation $(X^{-1}\t f^{-1})( Y^{-1} X^{-1} \t f)( X^{-1}Y^{-1} X^{-1} \t f^{-1})=1$. This is coherent with the calculation in \cite{L}.

\subsection{Comparison with \cite{FM}}

\begin{figure} 
\centerline{\relabelbox 
\epsfysize 11cm
\epsfbox{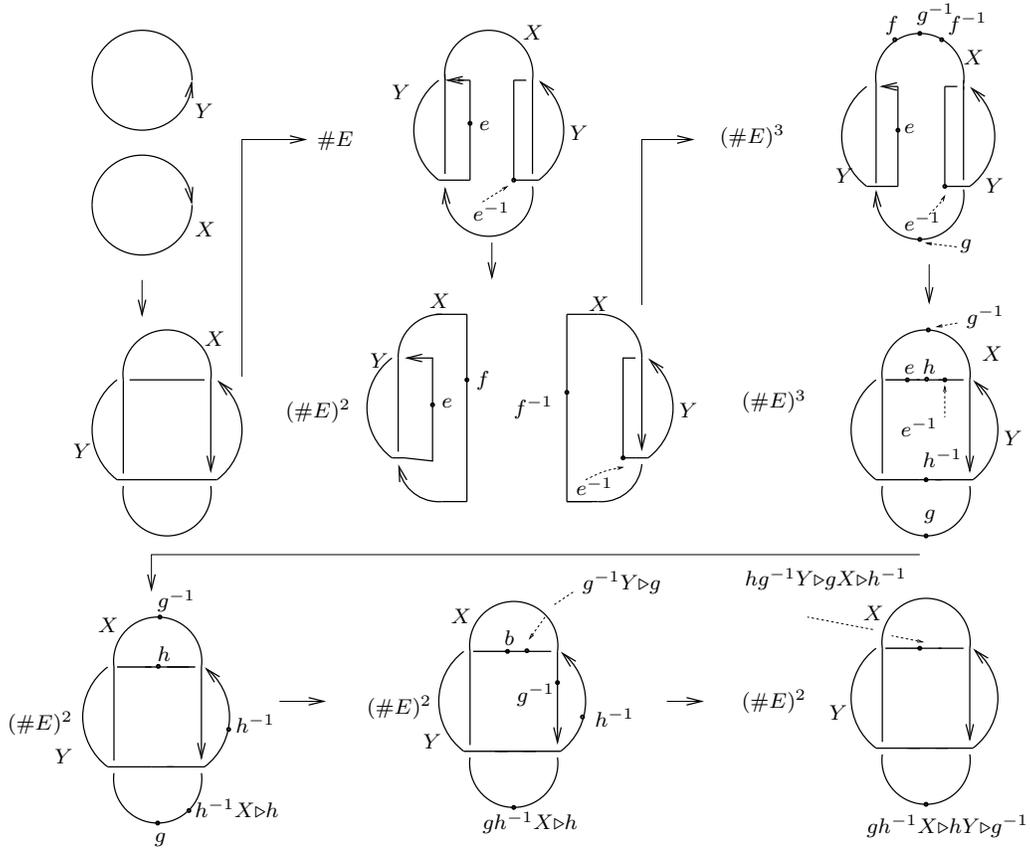}
\relabel {X10}{$\s{X}$}
\relabel {X}{$\s{X}$}
\relabel {X1}{$\s{X}$}
\relabel {X3}{$\s{X}$}
\relabel {X4}{$\s{X}$}
\relabel {X5}{$\s{X}$}
\relabel {X6}{$\s{X}$}
\relabel {X7}{$\s{X}$}
\relabel {X8}{$\s{X}$}
\relabel {X9}{$\s{X}$}
\relabel {Y}{$\s{Y}$}
\relabel {Y1}{$\s{Y}$}
\relabel {Y3}{$\s{Y}$}
\relabel {Y4}{$\s{Y}$}
\relabel {Y5}{$\s{Y}$}
\relabel {Y6}{$\s{Y}$}
\relabel {Y7}{$\s{Y}$}
\relabel {Y8}{$\s{Y}$}
\relabel {Y9}{$\s{Y}$}
\relabel {Y10}{$\s{Y}$}
\relabel {Y11}{$\s{Y}$}
\relabel {Y12}{$\s{Y}$}                  
\relabel {a}{$\s{e}$}
\relabel {b}{$\s{e^{-1}}$}
\relabel {g}{$\s{f}$}
\relabel {h}{$\s{f^{-1}}$}
\relabel {i}{$\s{e}$}
\relabel {j}{$\s{e^{-1}}$}
\relabel {k}{$\s{e}$}
\relabel {l}{$\s{e^{-1}}$}
\relabel {n}{$\s{f}$}
\relabel {p}{$\s{f^{-1}}$}
\relabel {m}{$\s{g}$}
\relabel {o}{$\s{g^{-1}}$}
\relabel {v}{$\s{g}$}
\relabel {t}{$\s{g^{-1}}$}
\relabel {u}{$\s{h}$}
\relabel {z}{$\s{h^{-1}}$}
\relabel {x}{$\s{h^{-1} X \t h}$}
\relabel {A}{$\s{b}$}
\relabel {B}{$\s{g^{-1}Y \t g }$}
\relabel {D}{$\s{g^{-1}}$}
\relabel {C}{$\s{gh^{-1} X \t h }$}
\relabel {E}{$\s{h^{-1}}$}
\relabel {F}{$\s{hg^{-1}Y \t g X \t h^{-1}}$}
\relabel {G}{$\s{g h^{-1}X \t h Y \t g^{-1}}$}
\relabel {M}{$\s{g^{-1}}$}
\relabel {N}{$\s{e}$}
\relabel {O}{$\s{h}$}
\relabel {P}{$\s{e^{-1}}$}
\relabel {Q}{$\s{h^{-1}}$}
\relabel {R}{$\s{g}$}
\relabel {s3}{$\s{\# E}$}
\relabel {s4}{$\s{(\# E)^2}$}
\relabel {s5}{$\s{(\# E)^3}$}
\relabel {s6}{$\s{(\# E)^3}$}
\relabel {s7}{$\s{(\# E)^2}$}
\relabel {s8}{$\s{(\# E)^2}$}
\relabel {s9}{$\s{(\# E)^2}$}
\endrelabelbox }
\caption{\label{calc2} Calculation of $J_\G(\S)$ where $\S$ is the Spun Hopf Link. We suppose that $G$ is abelian and that $\d=1_G$, which considerably simplifies the labellings. We also skip the summation signs and the very last step.}
\end{figure}

Let $\G=(G,E,\d,\t)$ be a finite crossed module and $J_\G$ the invariant of knotted surfaces defined in \cite{FM}. We can easily calculate $J_\G(\S)$, if $\S$ is the Spun Hopf Link, from figure \ref{movie}. We will make the calculation in the much simpler case $G$ and $E$ are abelian and $\d=1_G$. The calculation appears in figure \ref{calc2}. It then follows that: 
\begin{equation}\label{value2}
J_\G(\S)={\#\left \{X,Y \in G; h,g \in E: \substack{hg^{-1} (Y \t g)(X \t h^{-1})=1_E\\   \textrm{and}\\ gh^{-1} (X \t h)(Y \t g^{-1})=1_E}  \right\}}.
\end{equation}
Thus $J_\G(T)=I_\G(S^4 \setminus \n(\S))$ if $T$ is the Spun Hopf Link, and $\G=(G,E,\d,\t)$ is a crossed module with $G$ abelian and $\d=1_G$. To see this make the substitution $h \mapsto f^{-1}$ and $g \mapsto h$, and compare with equation (\ref{value1}).

Due to the calculations in this chapter, and the similarity between the two constructions, we make the following conjecture:

\begin{Conjecture}
Let $\S\subset S^4$ be an oriented knotted surface. Let also $\G=(G,E,\d,\t)$ be a finite crossed module. We have:

\begin{equation}
I_\G(S^4 \setminus \n(\S))=J_\G(\S)(\# E)^{-\chi (\S)},\end{equation}
where $J_\G$ is the invariant of knotted surfaces constructed in \cite{FM} and $\chi(\S)$ is the Euler Characteristic of $\S$.
\end{Conjecture}

\section*{Acknowledgements}
A large portion of this work was done while I was visiting the University of Nottingham in the spring of 2005. I have the  financial support of FCT (Portugal), post-doc grant number SFRH/BPD/17552/2004, part of the research project POCTI/MAT/60352/2004 (''Quantum Topology''), also financed by FCT.

I would like to thank  Roger Francis Picken, Marco Mackaay and Gustavo Granja for many helpful suggestions, and Scott Carter for encouraging comments, and for telling me about \cite{CKS}  and the approach there to handle decompositions of 2-knot complements.  I would also like to express my gratitude to Ronald Brown and Tim Porter for their generous help, suggestions and encouragement while I was visiting the University of Bangor, together with their hospitality; and to an anonymous referee of a previous version of this work, who, for example, suggested me the main  result  of corollary \ref{referee}.

\end{document}